\documentclass[10pt]{amsart}

\input{diagrams}

\newtheorem{proposition}{Proposition}[section]
\newtheorem{lemma}[proposition]{Lemma}
\newtheorem{corollary}[proposition]{Corollary}
\newtheorem{theorem}[proposition]{Theorem}

\theoremstyle{definition}

\newtheorem{example}[proposition]{Example}

\theoremstyle{remark}

\newcommand{\thlabel}[1]{\label{th:#1}}
\newcommand{\thref}[1]{Theorem~\ref{th:#1}}
\newcommand{\selabel}[1]{\label{se:#1}}
\newcommand{\seref}[1]{Section~\ref{se:#1}}
\newcommand{\lelabel}[1]{\label{le:#1}}
\newcommand{\leref}[1]{Lemma~\ref{le:#1}}
\newcommand{\prlabel}[1]{\label{pr:#1}}
\newcommand{\prref}[1]{Proposition~\ref{pr:#1}}
\newcommand{\colabel}[1]{\label{co:#1}}
\newcommand{\coref}[1]{Corollary~\ref{co:#1}}

\newcommand{\exlabel}[1]{\label{ex:#1}}

\newcommand{\eqlabel}[1]{\label{eq:#1}}
\newcommand{\equref}[1]{(\ref{eq:#1})}

\newcommand{\Az}{\dul{\rm Az}}
\newcommand{\coev}{{\rm coev}}
\newcommand{\ev}{{\rm ev}}
\newcommand{\Can}{{\rm Can}}
\newcommand{\Pic}{{\rm Pic}}
\newcommand{\Br}{{\rm Br}}
\newcommand{\Hom}{{\rm Hom}}
\newcommand{\End}{{\rm End}}

\newcommand{\Ker}{{\rm Ker}\,}

\newcommand{\im}{{\rm Im}\,}

\def\Ab{\underline{\underline{\rm Ab}}}
\def\lan{\langle}
\def\ran{\rangle}
\def\ot{\otimes}
\def\colim{{\rm colim}}

\def\RR{{\mathbb R}}

\def\CC{{\mathbb C}}

\def\units{{\mathbb G}_m}
\def\GG{{\mathbb G}}

\newcommand{\Cc}{\mathcal{C}}
\newcommand{\Dd}{\mathcal{D}}
\newcommand{\Ee}{\mathcal{E}}
\newcommand{\Ff}{\mathcal{F}}

\newcommand{\Mm}{\mathcal{M}}
\newcommand{\Pp}{\mathcal{P}}
\newcommand{\Ss}{\mathcal{S}}
\newcommand{\Rr}{\mathcal{R}}

\def\*C{{}^*\hspace*{-1pt}{\Cc}}
\def\text#1{{\rm {\rm #1}}}

\def\ol{\overline}
\def\ul{\underline}
\def\dul#1{\underline{\underline{#1}}}

\def\cat{\dul{\rm cat}}
\def\Ab{\dul{\rm Ab}}

\begin{document}
\title[The Brauer group of Azumaya corings ]{The Brauer group of Azumaya corings
and the second cohomology group}
\author{S. Caenepeel}
\address{Faculty of Engineering,
Vrije Universiteit Brussel, VUB, B-1050 Brussels, Belgium}
\email{scaenepe@vub.ac.be}
\urladdr{http://homepages.vub.ac.be/\~{}scaenepe/}
\author{B. Femi\'c}
\address{Faculty of Engineering,
Vrije Universiteit Brussel, VUB, B-1050 Brussels, Belgium}
\email{bfemic@vub.ac.be}
\thanks{}


\begin{abstract}
Let $R$ be a commutative ring. An Azumaya coring consists of a couple $(S,\Cc)$,
with $S$ a faithfully flat commutative $R$-algebra, and an $S$-coring $\Cc$
satisfying certain properties. If $S$ is faithfully projective, then the dual of $\Cc$
is an Azumaya algebra. Equivalence classes of Azumaya corings form an
abelian group, called the Brauer group of Azumaya corings. This group is
canonically isomorphic to the second flat cohomology group. We also give
algebraic interpretations of the second Amitsur cohomology group and the
first Villamayor-Zelinsky cohomology group in terms of corings.
\end{abstract}

\maketitle

\section*{Introduction}
Let $k$ be a field, and $l$ a Galois field extension of $k$ with group $G$. The Crossed
Product Theorem states that we have an isomorphism
$\Br(l/k)\cong H^2(G,l^*)$. The map from the second cohomology group to the
Brauer group can be described easily and explicitly: if $f\in Z ^2(G,l^*)$ is a 
$2$-cocycle, then the central simple algebra representing the class in
$\Br(l/k)$ corresponding to $f$ is
$$A=\bigoplus_{\sigma\in G} Au_{\sigma},$$
with multiplication rule
$$(au_\sigma)(bu_\tau)=a\sigma(b)f(\sigma,\tau)u_{\sigma\tau}.$$
From the fact that every central simple algebra can be split by a Galois
extension, it follows that the full Brauer group $\Br(k)$ can be
described as a second cohomology group
$$\Br(k)\cong H^2({\rm Gal}(k^{\rm sep}/k),k^{\rm sep*}),$$
where $k^{\rm sep}$ is the separable closure of $k$.\\
The definition of the Brauer group can be generalized from fields
to commutative rings (see \cite{AG}), or, more generally, to schemes
(see \cite{Gr2}). The cohomological description of the Brauer group of
a commutative ring is more complicated; first of all, Galois cohomology
is no longer sufficient, since not every Azumaya algebra can be split
by a Galois extension. More general cohomology theories have to be
introduced, such as Amitsur cohomology (over commutative rings) or
\v Cech cohomolgy (over schemes). The Crossed Product Theorem is
replaced by a long exact sequence, called the Chase-Rosenberg sequence.
We can introduce the second \'etale cohomology group $H^2(R_{\rm et},\units)$,
as the second right derived functor of a global section functor. If
$R=k$ is a field, then this group equals the total Galois cohomology group
$H^2({\rm Gal}(k^{\rm sep}/k),k^{\rm sep*})$. Then we have a monomorphism
$$\Br(R)\hookrightarrow H^2(R_{\rm et},\units).$$
In general, this monomorphism is not surjective, as the Brauer group is allways
torsion, and the second cohomology group is not torsion in general. Gabber
\cite{Gabber} proved that the Brauer group is isomorphic to the torsion part of the
second cohomology group.\\
In \cite{Ta}, Taylor introduced a new Brauer group, consisting of equivalence
classes of algebras that do not necessarily have a unit. The classical Brauer group
is a subgroup, and it is shown in \cite{RT} that Taylor's Brauer group is isomorphic to the full second
\'etale cohomology group. The proof depends on deep results, such as Artin's
Refinement Theorem (see \cite{A1}); also the proof does not provide an explicit
procedure producing a Taylor-Azumaya algebra out of an Amitsur cocycle.\\
In this paper, we propose a new Brauer group, and we show that it is isomorphic
to the full second flat cohomology group. The elements of this new Brauer
group are equivalence classes of corings. Corings were originally introduced by Sweedler
\cite{Sweedler75}; inspired by an observation made by Takeuchi that a large class
of generalized Hopf modules can be viewed as comodules over a coring,
Brzezi\'nski \cite{Brzezinski02} revived the theory of corings. 
\cite{Brzezinski02} was followed by a series of papers giving new applications of
corings, we refer to \cite{BrzezinskiWisbauer} for a survey.\\
Let $S$ be a commutative faithfully flat $R$-algebra. We can define a comultiplication
and a counit on the $S$-bimodule $S\ot S$, making $S\ot S$
into a coring. This coring, called Sweedler's canonical
coring, can be used to give an elegant approach to descent theory: the 
category of descent data is isomorphic to the category of comodules over the
coring. Our starting observation is now the following: an Amitsur $2$-cocycle can
be used to deform the comultiplication on $S\ot S$, such that the new comultiplication
is still coassociative. Thus the Amitsur $2$-cocycle condition should be viewed
as a coassociativity condition rather than an associativity condition (in contrast
with the Galois $2$-cocycle condition, which is really an associativity condition).
In the situation where $S$ is faithfully projective as an $R$-module, we can take
the dual of the coring $S\ot S$, which is an $S$-ring, isomorphic to $\End_R(S)$.
Amitsur $2$-cocycles can then be used to deform the multiplication on $\End_R(S)$,
leading to an Azumaya algebra in the classical sense; this construction leads to a map
$H^2(S/R,\units)\to \Br(S/R)$, and we will show that it is one of the maps in the Chase-Rosenberg
sequence. The duality between the $S$-coring $S\ot S$ and the $S$-ring
$\End_R(S)$ works well in both directions if $S/R$ is faithfully projective, but
fails otherwise; this provides an explanation for the fact that we need the condition that $S/R$
is faithfully projective in order to fit the relative Brauer group $\Br(S/R)$ into the
Chase-Rosenberg sequence.\\
The canonical coring construction can be generalized slightly: if $I$ is an invertible
$S$-module, then we can define a coring structure on $I^*\ot I$. Such a coring will
be called an elementary $S/R$-coring. Azumaya $S/R$-corings are then introduced
as twisted forms of elementary $S/R$-corings. If $S/R$ is faithfully projective, then
the dual of an Azumaya $S/R$-coring is an Azumaya algebra containing
$S$ as a maximal commutative subalgebra. The set of isomorphism classes of
Azumaya $S/R$-corings forms a group; after we divide by the subgroup consisting
of elementary corings, we obtain the relative Brauer group $\Br^c(S/R)$; we will
show that $\Br^c(S/R)$ is isomorphic to Villamayor and Zelinsky's cohomology
group with values in the category of invertible modules $H^1(S/R,\dul{\Pic})$
\cite{ViZ}. As a consequence, $\Br^c(S/R)$ fits into a Chase-Rosenberg type
sequence (even if $S/R$ is not faithfully projective).\\
An Azumaya coring will consist of a couple $(S,\Cc)$, where $S$ is a
(faithfully flat) commutative ring extension of $R$, and $\Cc$ is an $S/R$-coring.
On the set of isomorphism classes, we define a Brauer equivalence relation,
and show that the quotient set is a group under the operation induced by the
tensor product over $R$. This group is called the Brauer group of Azumaya corings,
and we can show that it is isomorphic to the full second cohomology group.\\
If $C$ is an object of a category $\Cc$, then the identity endomorphism of $C$
will also be denoted by $C$.

\section{The Brauer group of a commutative ring}\selabel{1}
\subsection{Amitsur cohomology}\selabel{1.1}
Let $R$ be a commutative ring, and $S$ an $R$-algebra that is faithfully flat
as an $R$-module. Tensor products over $R$ will be written without index $R$:
$M\ot N=M\ot_RN$, for $R$-modules $M$ and $N$. The $n$-fold tensor product
$S\ot\cdots \ot S$ will be denoted by $S^{\ot n}$. For $i\in \{1,\cdots,n+2\}$,
we have an algebra map
$$\eta_i:\ S^{\ot (n +1)}\to S^{\ot (n+2)},$$
given by
$$\eta_i(s_1\ot \cdots \ot s_{n+1})=s_1\ot \cdots\ot s_{i-1}\ot 1\ot s_i\ot\cdots\ot  s_{n+1}.$$
Let $P$ be a covariant functor from a full subcategory of the category of commutative
$R$-algebras that contains all tensor powers $S^{\ot n}$ of $S$ to abelian groups.
Then we consider
$$\delta_n=\sum_{i=1}^{n+2} (-1)^{i-1}P(\eta_i):\ P(S^{\ot (n +1)})\to
P(S^{\ot (n +2)}).$$
It is straightforward to show that $\delta_{n+1}\circ \delta_n=0$, so we obtain a complex
$$0\to P(S)\rTo^{\delta_0} P(S^{\ot 2})\rTo^{\delta_1} P(S^{\ot 3})\rTo^{\delta_2}\cdots,$$
called the Amitsur complex $\Cc(S/R)$. We write
$$Z^n(S/R, P)=\Ker\delta_n~~;~~B^n(S/R, P)=\im\delta_{n-1};$$
$$H^n(S/R, P)=Z^n(S/R, P)/B^n(S/R, P).$$
$H^n(S/R, P)$ will be called the $n$-th Amitsur cohomology group of $S/R$
with values in $P$. Elements in $Z^n(S/R,P)$ are called $n$-cocycles, and elements
in $B^n(S/R,P)$ are called $n$-coboundaries.\\
In this paper, we will mainly look at the following two examples: $P=\Pic$, where
$\Pic(S)$ is the Picard group of $S$, consisting of isomorphism classes of invertible
$S$-modules, and $P=\units$, where $\units(S)$ is the group consisting of all
invertible elements of $S$.\\
If $u\in S^{\ot n}$, then we will write $u_i=\eta_i(u)$. Observe that $u\in
\units(S^{\ot 3})$ is then a cocycle in $Z^2(S/R,\units)$ if and only if
$$u_1u_2^{-1}u_3u_4^{-1}=1.$$
Amitsur cohomology was first introduced in \cite{Amitsur1} (over fields); it can be viewed as an
affine version of \v Cech cohomology. For a more detailed discussion, see for example
\cite{Caenepeel98,CR,KO1}. We now present some elementary properties of
Amitsur cohomology groups. We will adopt the following notation: an element $u\in S^{\ot n}$
will be written formally as
$u=u^1\ot u^2\ot \cdots \ot u^n$, 
where the summation is understood implicitly.

\begin{proposition}\prlabel{am1a}
Let $R$ be a commutative ring, and $f:\ S\to T$ a morphism of commutative
$R$-algebras. $f$ induces maps $f_*:\ H^n(S/R,P)\to H^n(T/R,P)$.
If $g:\ S\to T$ is a second algebra map, then $f_*=g_*$ (for $n\geq 1$).
\end{proposition}

\begin{proof}
The first statement is obvious. For the proof of the second one, we refer to
\cite[Prop. 5.1.7]{KO1}.
\end{proof}

The following result is obvious.

\begin{lemma}\lelabel{am1}
If $u,v\in Z^n(S/R,\units)$, then
$$u\ot v= (u^1\ot v^1)\ot (u^2\ot v^2)\ot\cdots\ot (u^n\ot v^n)\in Z^n(S\ot S/R,\units).$$
If $u,v\in B^n(S/R,\units)$, then $u\ot v\in B^n(S\ot S/R,\units)$.
\end{lemma}

\begin{corollary}\colabel{am2}
If $u\in Z^n(S/R,\units)$, then
 $[u\ot 1]=[1\ot u]$, and
 $[u\ot u^{-1}]=1$ in $H^n(S\ot S/R,\units)$.
\end{corollary}

\begin{proof}
Apply \prref{am1a} to the algebra maps $\eta_1,\eta_2:\
S\to S\ot_R S$, $\eta_1(s)=1\ot s$, $\eta_2(s)=s\ot 1$.
\end{proof}

\begin{lemma}\lelabel{am3}
Take a cocycle $u=u^1\ot u^2\ot u^3\in Z^2(S/R,\units)$. $|u|=u^1u^2u^3\in \units(S)$
is called the norm of $u$, and
$$u_1\ot |u|^{-1}u_2u_3=1\ot 1=|u|^{-1}u_1u_2\ot u_3.$$
\end{lemma}

\begin{proof}
The first equality is obtained after we multiply the second, third and fourth tensor factors
in the cocycle condition $u_1u_2^{-1}u_3u_4^{-1}=1$. The second equality is obtained
after multiplying the first three tensor factors.
\end{proof}

A 2-cocycle $u$ is called normalized if $|u|=1$.

\begin{lemma}\lelabel{am4}
Every cocycle $u$ is cohomologous to a normalized cocycle.
\end{lemma}

\begin{proof}
First observe that $\Delta_1(|u|^{-1}\ot 1)=1\ot |u|^{-1}\ot 1$. The cocycle
$u\Delta_1(|u|^{-1}\ot 1)=u^1\ot |u|^{-1}u_2\ot u_3$ is normalized and cohomologous
to $u$.
\end{proof}

Now we consider the Amitsur complex $\Cc(S\ot S/R\ot S)$. We have a natural
isomorphism
$$(S\ot S)^{\ot_{R\ot S}n}\rTo^{\cong} S^{\ot (n+1)},~~
(s_1\ot t_1)\ot\cdots\ot (s_n\ot t_n)\mapsto s_1\ot\cdots\ot s_n\ot t_1\cdots t_n.$$
The augmentation maps ($i=1,2,3$)
$$\eta_i:\ (S\ot S)^{\ot_{R\ot S}2}\to (S\ot S)^{\ot_{R\ot S}3}$$
can then be viewed as maps
$$\eta_i:\ S^{\ot 3}\to S^{\ot 4},$$
and we find, for $u\in Z^2(S/R,\units)$ and $i=1,2,3$ that $\eta_i(u)=u_i$.
Consequently $u\ot 1=u_4=u_1u_2^{-1}u_3=\Delta_1(u)\in B^2(S\ot S/R\ot S,\units)$.

\begin{lemma}\lelabel{am5}
If $u\in Z^2(S/R,\units)$, then $u\ot 1\in B^2(S\ot S/R\ot S,\units)$.
\end{lemma}

\subsection{Derived functor cohomology}\selabel{1.2}
Let $R$ be a commutative ring. $\cat(R_{\rm fl})$ is the full subcategory of commutative
flat finitely presented $R$-algebras.
A covariant functor $P:\ \cat(R_{\rm fl})\to \Ab$ is called a presheaf on $R_{\rm fl}$.
The category of presheaves on $R_{\rm fl}$ and natural transformations will
be denoted by $\Pp(R_{\rm fl})$.\\
A presheaf $P$ is called a sheaf if $H^0(S'/S,P)=P(S)$, for every faithfully flat
$R$-algebra homomorphism $S\to S'$. The full subcategory of $\Pp(R_{\rm fl})$
consisting of sheaves is denoted by $\Ss(R_{\rm fl})$.
$\Pp(R_{\rm fl})$ and $\Ss(R_{\rm fl})$ are abelian categories having enough
injective objects.\\
$\GG_a$ and $\units$ are sheaves on $R_{\rm fl}$. The embedding functor
$i:\ \Ss(R_{\rm fl})\to \Pp(R_{\rm fl})$ has a left adjoint $a:\ \Pp(R_{\rm fl})\to \Ss(R_{\rm fl})$.\\
The ``global section" functor $\Gamma:\ \Ss(R_{\rm fl})\to \Ab$ is left exact, so
we can consider its $n$-th right derived functor $R^n\Gamma$. We define the
$n$-th flat cohomology group by
$$H^n(R_{\rm fl},\units)=R^n\Gamma(\units).$$
Fix a faithfully flat $R$-algebra $S$, and consider the functor
$$g=H^0(S/R,-):\ \Pp(R_{\rm fl})\to \Ab.$$
Then $\Gamma=g\circ i$, and $i$ takes injective objects of $\Ss(R_{\rm fl})$
to $g$-acyclics (see \cite[lemma 5.6.6]{Caenepeel98}), and we have long exact
sequences, for every sheaf $F$, and for every $q\geq 0$ (see \cite{Caenepeel98,ViZ}):
\begin{eqnarray}\eqlabel{et1}
0&\longrightarrow& H^1(S/R,C^q)\longrightarrow H^{q+1}(R_{\rm fl},F)
\longrightarrow H^0(S/R,H^{q+1}(\bullet,F))\\
&\longrightarrow& H^2(S/R,C^q)\longrightarrow H^1(S/R,C^{q+1})
\longrightarrow H^1(S/R,H^{q+1}(\bullet,F))\nonumber\\
&\longrightarrow& \cdots\nonumber\\
&\longrightarrow& H^{p+1}(S/R,C^q)\longrightarrow H^p(S/R,C^{q+1})
\longrightarrow H^p(S/R,H^{q+1}(\bullet,F))\nonumber\\
&\longrightarrow& \cdots\nonumber.
\end{eqnarray}
The sheaf $C^i$ is the $i$-th syzygy of an injective resolution
$0\to F\to X^0\to X^1\to \cdots$ of $F$ in $\Ss(R_{\rm fl})$, that is,
$C^i=\Ker(X^i\to X^{i+1})$.\\
A morphism $f:\ S\to T$ of commutative faithfully flat $R$-algebras
induces a map between the corresponding sequences \equref{et1}, namely
we have a commutative diagram
\begin{equation}\eqlabel{et1a}
\begin{diagram}
0\rightarrow&H^1(S/R,C^q)&\to&H^{q+1}(R_{\rm fl},F)&\to&
H^0(S/R,H^{q+1}(\bullet,F))&\to&\cdots \\
&\dTo_{f_*}&&\dTo_{=}&&\dTo_{f_*}&&\\
0\rightarrow&H^1(T/R,C^q)&\to&H^{q+1}(R_{\rm fl},F)&\to&
H^0(T/R,H^{q+1}(\bullet,F))&\to&\cdots
\end{diagram}
\end{equation}
It is known that $H^1(R_{\rm fl},\units)=\Pic(R)$, the group of rank one
projective $R$-modules. Writing down \equref{et1} for $F=\units$ and
$q=0$, we find the exact sequence
\begin{eqnarray}\eqlabel{et2}
0&\longrightarrow& H^1(S/R,\units)\longrightarrow \Pic(R)
\longrightarrow H^0(S/R,\Pic)\\
&\longrightarrow& H^2(S/R,\units)\longrightarrow H^1(S/R,C^{1})
\longrightarrow H^1(S/R,\Pic)\nonumber\\
&\longrightarrow&H^3(S/R,\units)\longrightarrow\cdots\nonumber
\end{eqnarray}
Let $\Rr$ be the category with faithfully flat commutative $R$-algebras as
objects. The set of morphisms between two objects $S$ and $T$ is a singleton
if there exists an algebra morphism $S\to T$ (then we write $S\leq T$), and
is empty otherwise. Then $\Rr$ is a directed preorder, that is a category
with at most one morphism between two objects, and such that every pair
of objects $(S,T)$ has a successor, namely $S\ot T$ (see \cite[IX.1]{McLane}).\\
Let $P$ be a presheaf on $R_{\rm fl}$. It follows from \prref{am1a}
that we have a functor
$$H^n(\bullet/R,P):\ \Rr\to \Ab,$$
and we can consider the colimit
$$\check H^n(R_{\rm fl},P)=\colim H^n(\bullet/R,P).$$
Now let $F$ be a sheaf. Using the exact sequences \equref{et1} and the
commutative diagrams \equref{et1a}, we find a homomorphism of abelian groups
$$\check H^n(R_{\rm fl},F)\to H^n(R_{\rm fl},F).$$
If $n=1$, this map is an isomorphism. In particular, we have
\begin{equation}\eqlabel{et3}
H^2(R_{\rm fl},\units)\cong H^1(R_{\rm fl},C^1)\cong\check H^1(R_{\rm fl},C^1).
\end{equation}
The category $\cat(R_{\rm fl})$ can be replaced by $\cat(R_{\rm et})$, the category
of \'etale $R$-algebras. All results remain valid,
and, moreover, we have
$$\check H^n(R_{\rm et},F)\cong H^n(R_{\rm et},F).$$
The proof of this result is based on Artin's Refinement Theorem \cite{A1}.

\subsection{Amitsur cohomology with values in $\dul{\Pic}$}\selabel{1.2b}
Let $R$ be a commutative ring. The category of invertible $R$-modules and
$R$-module isomorphisms is denoted by $\dul{\Pic}(R)$. The Grothendieck group
$K_0\dul{\Pic}(R)$ is the Picard group $\Pic(R)$.
The inverse of $[I]\in \Pic(R)$ is represented by $I^*=\Hom_R(I,R)$.
If $I\in \dul{\Pic}(R)$, then the evaluation map
$\ev_I:\ I\ot I^*\to R$ is an isomorphism, with inverse the coevaluation map
$\coev_I:\ R\to I\ot I^*$. If $\coev_I(1)=\sum_i e_i\ot e_i^*$, then
$\{(e_i,e_i^*)~|~i=1,\cdots n\}$ is a finite dual basis for $I$.\\ 

Let $S$ be a commutative faithfully flat $R$-algebra. For every positive integer
$n$, we have a functor
$$\delta_{n-1}:\ \dul{\Pic}(S^{\ot n})\to \dul{\Pic}(S^{\ot (n+1)}),$$
given by
$$\delta_{n-1}(I)=I_1\ot_{S^{\ot (n+1)}}I^*_2\ot_{S^{\ot (n+1)}}\cdots
\ot_{S^{\ot (n+1)}}J_{n+1},$$
$$\delta_{n-1}(f)=f_1\ot_{S^{\ot (n+1)}} (f^*_2)^{-1}\ot_{S^{\ot (n+1)}}\cdots
\ot_{S^{\ot (n+1)}} (g_{n+1})^{\pm 1},$$
with $J=I$ or $I^*$, $g=f$ or $f^*$ depending on whether $n$ is odd or even.
Here $I_i= I\ot_{S^{\ot n}}S^{\ot {n+1}}$, where $S^{\ot {n+1}}$ is a left
$S^{\ot n}$-module via $\eta_i:\ S^{\ot n}\to S^{\ot {n+1}}$ (see \seref{1.1}). 
We easily compute that
$$\delta_n\delta_{n-1}(I)=\bigotimes_{j=1}^{n+2}\bigotimes_{i=1}^{j-1}
(I_{ij}\ot_{S^{\ot (n+2)}} I^*_{ij}),$$
so we have a natural isomorphism
$$\lambda_I= \bigotimes_{j=1}^{n+2}\bigotimes_{i=1}^{j-1} {\rm ev}_{I_{ij}}:\
\delta_n\delta_{n-1}(I)\to S^{\ot (n+2)}.$$
$\dul{Z}^{n-1}(S/R,\dul{\Pic})$ is the category with objects $(I,\alpha)$,
with $I\in \dul{\Pic}(S^{\ot n})$, and $\alpha:\ \delta_{n-1}(I)\to S^{\ot(n+1)}$
an isomorphism of $S^{\ot(n+1)}$-modules such that $\delta_n(\alpha)=\lambda_I$.
A morphism $(I,\alpha)\to (J,\beta)$ is an isomorphism of $S^{\ot n}$-modules
$f:\ I\to J$ such that $\beta\circ \delta_{n-1}(f)=\alpha$. 
$\dul{Z}^{n-1}(S/R,\dul{\Pic})$ is a symmetric monoidal category, with tensor
product $(I,\alpha)\ot (J,\beta)=(I\ot_{S^{\ot n}}J,\alpha\ot_{S^{\ot (n+1)}}\beta)$
and unit object $(S^{\ot n},S^{\ot (n+1)})$. Every object in this category is
invertible, and we can consider
$$K_0\dul{Z}^{n-1}(S/R,\dul{\Pic})={Z}^{n-1}(S/R,\dul{\Pic}).$$
We have a strongly monoidal functor
$$\delta_{n-2}:\ \dul{\Pic}(S^{\ot(n-1)})\to \dul{Z}^{n-1}(S/R,\dul{\Pic}),$$
$\delta_{n-2}(J)=(\delta_{n-2}(J),\lambda_J)$. Consider the
subgroup $B^{n-1}(S/R,\dul{\Pic})$ of  $Z^{n-1}(S/R,\dul{\Pic})$,
consisting of elements represented by $\delta_{n-2}(J)$, with $J\in\dul{\Pic}(S^{\ot n-1})$.
We then define
$$ H^{n-1}(S/R,\dul{\Pic})=Z^{n-1}(S/R,\dul{\Pic})/B^{n-1}(S/R,\dul{\Pic}).$$
This definition is such that we have a long exact sequence
(see \cite{ViZ}):
\begin{eqnarray}\eqlabel{villa1}
0&\longrightarrow& H^1(S/R,\units)\longrightarrow \Pic(R)
\longrightarrow H^0(S/R,\Pic)\\
&\longrightarrow& H^2(S/R,\units)\longrightarrow H^1(S/R,\dul{\Pic})
\longrightarrow H^1(S/R,\Pic)\nonumber\\
&\longrightarrow& \cdots\nonumber\\
&\longrightarrow& H^{p+1}(S/R,\units)\longrightarrow H^p(S/R,\dul{\Pic})
\longrightarrow H^p(S/R,\Pic)\nonumber\\
&\longrightarrow& \cdots\nonumber.
\end{eqnarray}
Comparing to \equref{et1} in the situation where $F=\units$ and $q=0$, we see that
\begin{equation}\eqlabel{villa2}
H^n(S/R,\dul{\Pic})\cong H^n(S/R,C^1),
\end{equation}
for all $n\geq 1$. For detail, we refer to \cite{Caenepeel98,ViZ}. The following
result can be viewed as an analog of \leref{am5}.

\begin{lemma}\lelabel{1.2a.2}
Let $(I,\alpha)\in\dul{Z}^1(S/R,\dul{\Pic})$. Then
$$(I\ot S,\alpha\ot S)\cong \delta_0(I)~~{\rm in}~~\dul{Z}^1(S\ot S/R\ot S,\dul{\Pic}),$$
and consequently $[I\ot S,\alpha\ot S]=1$ in $H^1(S\ot S/R\ot S,\dul{\Pic})$.
\end{lemma}

\begin{proof}
The isomorphism $\alpha:\ I_1\ot_{S^{\ot 3}}I^*_2\ot_{S^{\ot 3}}I_3\to S^{\ot 3}$
induces an isomorphism
$$\beta:\ I_3=I\ot S\to I^*_1\ot_{S^{\ot 3}}I_2=
(S\ot I)^*\ot_{(S\ot S)\ot_{R\ot S}(S\ot S)}(S\ot I).$$
The fact that $\delta_2(\alpha)=\lambda_I$ implies that $\beta$ is an 
isomorphism in $\dul{Z}^1(S\ot S/R\ot S,\dul{\Pic})$.
\end{proof}

\begin{proposition}\prlabel{1.2a.3}
Let $f:\ S\to T$ be a morphism of commutative faithfully flat $R$-algebras.
$f$ induces group morphisms $f_*:\ H^n(S/R,\dul{\Pic})\to
H^n(T/R,\dul{\Pic})$. If $g:\ S\to T$ is a second algebra morphism, then
$f_*=g_*$.
\end{proposition}

\begin{proof}
We have a functor $f_*:\ \dul{Z}^{n-1}(S/R,\dul{\Pic})\to
\dul{Z}^{n-1}(T/R,\dul{\Pic})$, given by
$$f_*(I,\alpha)=(I\ot_{S^{\ot n}}T^{\ot n},\alpha\ot_{S^{\ot n+1}}T^{\ot n+1}).$$
$f_*$ induces maps $f_*:\ H^n(S/R,\dul{\Pic})\to
H^n(T/R,\dul{\Pic})$.\\
$f$ and $g$ induce maps $f_*$ and $g_*$
between the exact sequence \equref{villa1} and its analog with $S$ replaced
by $T$. We have seen in \prref{am1a} that these maps coincide on
$H^n(S/R,\units)$ and $H^n(S/R,\Pic)$. It follows from the five lemma that
they also coincide on $H^n(T/R,\dul{\Pic})$.
\end{proof}

It follows from \prref{1.2a.3} that we have a functor
$$H^1(\bullet/R,\dul{\Pic}):\ \Rr\to \Ab,$$
so we can consider the colimit
$$\check H^n(R_{\rm fl},\dul{\Pic})=\colim H^1(\bullet/R,\dul{\Pic}).$$
If $f:\ S\to T$ is a morphism of commutative faithfully flat $R$-algebras,
then the maps $f_*$ establish a map between the corresponding exact sequences
\equref{villa1}. This implies that the isomorphisms \equref{villa2} fit into
commutative diagrams
$$\begin{diagram}
H^n(S/R,\dul{\Pic})&\rTo^{\cong}&H^n(S/R,C^1)\\
\dTo^{f_*}&&\dTo_{f_*}\\
H^n(T/R,\dul{\Pic})&\rTo^{\cong}&H^n(T/R,C^1)
\end{diagram}$$
Consequently, the functors $H^n(\bullet/R,\dul{\Pic})$ and $H^n(\bullet/R,C^1)$
are isomorphic, and
\begin{equation}\eqlabel{villa10}
\check H^1(R_{\rm fl},\dul{\Pic})\cong \check H^1(R_{\rm fl},C^1)
\cong H^2(R_{\rm fl},\units).
\end{equation}

\subsection{The Brauer group}\selabel{1.3}
Let $R$ be a commutative ring. An $R$-algebra $A$ is called an Azumaya
algebra if there exists a commutative faithfully flat $R$-algebra $S$ such that $A\ot S
\cong \End_S(P)$ for some faithfully projective $S$-module $P$. There are
several equivalent characterizations of Azumaya algebras, we refer to the literature
\cite{Caenepeel98,DI,KO1}. An Azumaya algebra over a field is nothing else then
a central simple algebra.\\
Two $R$-Azumaya algebras $A$ and $B$ are called Brauer equivalent if there exist
faithfully projective $R$-modules $P$ and $Q$ such that $A\ot \End(P)\cong
B\ot \End(Q)$ as $R$-algebras. This induces an equivalence relation on the set of
isomorphism classes of $R$-Azumaya algebras. The quotient set $\Br(R)$ is
an abelian group under the operation induced by the tensor product. The inverse
of a class represented by an algebra $A$ is represented by the opposite algebra
$A^{\rm op}$.\\
If $i:\ R\to S$ is a morphism of commutative rings, then we have an associated
abelian group map 
$$\Br(i):\ \Br(R)\to \Br(S),~~i[A]=[A\ot S].$$
The kernel $\Ker(\Br(i))=\Br(S/R)$ is called the part of the Brauer group of $R$
split by $S$.\\
If $S/R$ is faithfully flat, then we have an embedding $\Br(S/R)\to H^1(S/R,C^1)$.
This embedding is an isomorphism if $S$ is faithfully projective as an $R$-module.
Consequently, we have an embedding
$$\Br(S/R)\to H^2(R_{\rm fl},\units),$$
and
$$\Br(R)\to H^2(R_{\rm fl},\units).$$
Since every $R$-Azumaya algebra can be split by an \'etale covering, $H^2(R_{\rm fl},\units)$
can be replaced by $H^2(R_{\rm et},\units)$ in the two formulas above. If $R$
is a field, or, more generally, if $R$ is a regular ring, then we have an
isomorphism
$$\Br(R)\cong H^2(R_{\rm et},\units).$$
In general, we do not have such an isomorphism, because the Brauer group is
torsion, and the second cohomology group is not (see \cite{Gr2}). Gabber
(\cite{Gabber}, see also \cite{KO3}) showed that
$$\Br(R)\cong H^2(R_{\rm et},\units)_{\rm tors},$$
for every commutative ring $R$. Taylor \cite{Ta} introduced a Brauer group
$\Br'(R)$ consisting of classes of algebras that have not necessarily a unit, but
satisfy a weaker property. $\Br'(R)$ contains $\Br(R)$ as a subgroup, and we have
an isomorphism \cite{RT}
$$\Br'(R)\cong H^2(R_{\rm et},\units).$$
The proof is technical, and relies on Artin's refinement Theorem \cite{A1}. It provides
no explicit description of the Taylor-Azumaya algebra that corresponds to a given
cocycle.

\section{Some adjointness properties}\selabel{2.1}
We start this technical Section with the following elementary fact.
For any morphism $\eta:\ R\to S$ of rings, we have an adjoint pair of functors
$(F=-\ot_R S,G)$ between the module categories $\Mm_R$ and $\Mm_S$.
$F$ is called the induction functor, and $G$ is the restriction of scalars functor.
For every $M\in \Mm_R$, $N\in \Mm_S$, we have a natural isomorphism
$$\Hom_R(M,G(N))\cong \Hom_S(M\ot_R S,N).$$
$f:\ M\to G(N)$ and the corresponding $\tilde{f}:\ M\ot_RS\to N$ are related by
the following formula:
\begin{equation}\eqlabel{2.1.0.1}
\tilde{f}(m\ot_R s)=f(m)s.
\end{equation}

Now assume that $R$ and $S$ are commutative rings, and consider the
ring morphisms $\eta_i:\ S\ot_R S\to S\ot_RS\ot_RS$ ($i=1,2,3$) introduced at the beginning
of \seref{1.1}. The corresponding adjoint pairs of functors between
$\Mm_{S^{\ot 2}}$ and $\Mm_{S^{\ot 3}}$ will be written as $(F_i,G_i)$.
$M\in \Mm_{S^{\ot 2}}$ will also be regarded as an $S$-bimodule, and we will
denote $M_i=F_i(M)$. For $m\in M$, we write
$$m_i=(M\ot_{S^{\ot 2}}\eta_i)(m).$$
In particular, $m_3=m\ot 1$ and $m_1=1\ot m$.

\begin{lemma}\lelabel{2.1}
Let $M\in \Mm_{S^{\ot 2}}$. Then we have an $S$-bimodule isomorphism
$$G_2(M_3\ot_{S^{\ot 3}}M_1)\cong M\ot_S M,$$
and an isomorphism
$${}_S\Hom_S(M,M\ot_S M)\cong \Hom_{S^{\ot 3}}(M_2, M_3\ot_{S^{\ot 3}}M_1).$$
\end{lemma}

\begin{proof}
The map
$$\alpha:\ M_3\ot M_1\to M\ot_S M,~~
\alpha((m\ot s)\ot (t\ot n))=tm\ot_S ns$$
induces a well-defined map
$$\alpha:\ M_3\ot_{S^{\ot 3}}M_1\to M\ot_S M.$$
Indeed, for all $m,n\in M$ and $s,t,u,v,w\in S$, we easily compute that
\begin{eqnarray*}
&&\hspace*{-2cm}
\alpha\bigl((m\ot s)(u\ot v\ot w)\ot (t\ot n)\bigr)
=\alpha\bigl((umv\ot sw)\ot (t\ot n)\bigr)\\
&=& tumv\ot_S nsw=utm\ot_S vnws\\
&=& \alpha\bigl((m\ot s)\ot (ut\ot vnw)\bigr)
=\alpha\bigl((m\ot s)\ot (u\ot v\ot w) (t\ot n)\bigr).
\end{eqnarray*}
The map
$$\beta:\ M\ot M\to M_3\ot_{S^{\ot 3}}M_1,~~\beta(m\ot n)=m_3\ot_{S^{\ot 3}} m_1$$
induces a well-defined map
$$\beta:\ M\ot_S M\to M_3\ot_{S^{\ot 3}}M_1.$$
Indeed,
\begin{eqnarray*}
&&\hspace*{-2cm}
\beta(ms\ot n)=(ms\ot 1)\ot_{S^{\ot 3}}(1\ot n)\\
&=&
(m\ot 1)(1\ot s\ot 1)\ot_{S^{\ot 3}}(1\ot n)\\
&=& (m\ot 1)\ot_{S^{\ot 3}}(1\ot s\ot 1)(1\ot n)\\
 &=&(m\ot 1)\ot_{S^{\ot 3}}(1\ot sn)
=\beta(m\ot sn).
\end{eqnarray*}
It is clear that $\alpha$ and $\beta$ are inverse $S$-bimodule maps. Finally,
the adjunction cited above tells us that
$${}_S\Hom_S(M,M\ot_S M)\cong \Hom_{S^{\ot 2}}(M, G_2(M_3\ot_{S^{\ot 3}}M_1))
\cong\Hom_{S^{\ot 3}}(M_2, M_3\ot_{S^{\ot 3}}M_1).$$
\end{proof}

Using \equref{2.1.0.1}, we can write an explicit formula for the map
$\tilde{f}:\ M_2\to M_3\ot_{S^{\ot 3}}M_1$ corresponding to
$f:\ M\to M\ot_S M$. To this end, we first introduce the following Sweedler-type
notation:
$$f(m)=m_{(1)}\ot_S m_{(2)},$$
where summation is understood implicitly. Then we have
\begin{equation}\eqlabel{2.1.1}
\tilde{f}(m_2)=\beta(f(m))=m_{(1)3}\ot_{S^{\ot 3}} m_{(2)1}.
\end{equation}

For $i=1,2,3,4$ and $j=1,2,3$, we now consider the ring morphisms
$$\eta_{ij}=\eta_i\circ \eta_j:\ S\ot_RS\to S\ot_RS\ot_RS\ot_RS$$
and the corresponding pairs of adjunct functors $(F_{ij},G_{ij})$ between the categories
$\Mm_{S^{\ot 2}}$ and $\Mm_{S^{\ot 4}}$.

\begin{lemma}\lelabel{2.1a}
Let $M\in \Mm_{S^{\ot 2}}$. Then we have a natural isomorphism of $S$-bimodules
$$G_{23}(M_{34}\ot_{S^{\ot 4}}M_{14}\ot_{S^{\ot 4}}M_{12})\cong
M\ot_SM\ot_SM,$$
and an isomorphism
$${}_S\Hom_S(M,M\ot_S M\ot_S M)\cong 
\Hom_{S^{\ot 4}}(M_{23}, M_{34}\ot_{S^{\ot 4}}M_{14}\ot_{S^{\ot 4}}M_{12}).$$
The map $\tilde{f}$ corresponding to $f\in {}_S\Hom_S(M,M\ot_S M\ot_S M)$, with
$f(m)= m_{(1)}\ot_S m_{(2)}\ot_S m_{(3)}$ is given by the formula
\begin{equation}\eqlabel{2.1a.1}
\tilde{f}(m_{23})=m_{(1)34}\ot_{S^{\ot 4}} m_{(2)14}\ot_{S^{\ot 4}} m_{(3)12}.
\end{equation}
\end{lemma}

\begin{proof}
The map
$$\alpha:\ M_{34}\ot_{S^{\ot 4}}M_{14}\ot_{S^{\ot 4}}M_{12}\to 
M\ot_S M\ot_S M$$
and
$$\beta:\ M\ot_S M\ot_S M\to M_{34}\ot_{S^{\ot 4}}M_{14}\ot_{S^{\ot 4}}M_{12}$$
given by the formulas
$$\alpha\Bigl((m\ot s\ot t)\ot_{S^{\ot 4}}(s'\ot n\ot t')\ot_{S^{\ot 4}}(s''\ot t''\ot p)\Bigr)=
s''s'm\ot_S t''ns\ot_S ptt'$$
and
$$\beta(m\ot_S n\ot_S p)=m_{34}\ot_{S^{\ot 4}} n_{14}\ot_{S^{\ot 4}}p_{12}$$
are well-defined inverse $S$-bimodule maps. Verification of the details goes
precisely as in the proof of \leref{2.1}. Then, using the adjunction from the
beginning of this Section, we find
\begin{eqnarray*}
&&\hspace*{-2cm}
{}_S\Hom_S(M,M\ot_S M\ot_S M)\cong \Hom_{S^{\ot 2}}(M, 
G_{23}(M_{34}\ot_{S^{\ot 4}}M_{14}\ot_{S^{\ot 4}}M_{12})\\
&\cong &
\Hom_{S^{\ot 4}}(M_{23}, M_{34}\ot_{S^{\ot 4}}M_{14}\ot_{S^{\ot 4}}M_{12}).
\end{eqnarray*}
Using \equref{2.1.0.1}, we find that $\tilde{f}(m_{23})=\beta(f(m))$, and
\equref{2.1a.1} then follows easily.
\end{proof}

Let $S$ be a commutative faithfully flat $R$-algebra. We have an algebra morphism
$m:\ S^{\ot n}\to S$, $m(s_1\ot\cdots\ot s_n)=s_1\cdots s_n$, and the
corresponding induction functor
$$-\ot_{S^{\ot n}} S=|-|:\ \Mm_{S^{\ot n}}\to \Mm_S,$$
which is strongly monoidal since $|S^{\ot n}|=S$, and
\begin{eqnarray*}
&&\hspace{-2cm} |M\ot_{S^{\ot n}}N|= M\ot_{S^{\ot n}}N\ot_{S^{\ot n}}S\\
&\cong& (M\ot_{S^{\ot n}} S)\ot_S (N\ot_{S^{\ot n}} S)\cong |M|\ot_S|N|.
\end{eqnarray*}
Recall from \cite[IX.4.6]{Ba2} that an $R$-module $M$ is faithfully projective if and
only if there exists an $R$-module $N$ such that $M\ot N\cong R^m$. This implies
that $|-|$ sends faithfully projective (resp. invertible) $S^{\ot n}$-modules to faithfully projective (resp. invertible) $S$-modules.

\begin{lemma}\lelabel{2.2a}
Let $M_1,\cdots, M_n \in \Mm_S$. Then
$$|M_1\ot\cdots\ot M_n|\cong M_1\ot_S\cdots\ot_S M_n.$$
\end{lemma}

\begin{proof}
The natural epimorphism $\pi:\ M_1\ot\cdots\ot M_n\to |M_1\ot\cdots\ot M_n|$
factors through $M_1\ot_S\cdots\ot_S M_n$ since
\begin{eqnarray*}
&&\hspace*{-20mm}
\pi(m_1\ot\cdots \ot sm_i\ot \cdots\ot m_n)=
(m_1\ot\cdots \ot sm_i\ot \cdots\ot m_n)\ot_{S^{\ot n}} 1\\
&=& (m_1\ot\cdots \ot m_i\ot \cdots\ot m_n)\ot_{S^{\ot n}} s\\
&=&(m_1\ot\cdots \ot sm_j\ot \cdots\ot m_n)\ot_{S^{\ot n}} 1,
\end{eqnarray*}
for all $i,j$, so we have a map
$$\alpha:\ M_1\ot_S\cdots\ot_S M_n\to |M_1\ot\cdots\ot M_n|.$$
In a similar way, the quotient map $M_1\ot\cdots\ot  M_n\to M_1\ot_S\cdots\ot_S M_n$
factors through $|M_1\ot\cdots\ot M_n|$, so we have a map
$$\beta:\ |M_1\ot\cdots\ot  M_n|\to M_1\ot_S\cdots\ot_S M_n,$$
which is inverse to $\alpha$.
\end{proof}

\section{Corings}\selabel{2.2}
Let $S$ be a ring. Recall that an $S$-coring is a coalgebra (or comonoid) $\Cc$
in the category ${}_S\Mm_S$. This means that $\Cc$ is an $S$-bimodule,
together with two $S$-bimodule maps $\Delta_\Cc:\ \Cc\to \Cc\ot_S \Cc$ and
$\varepsilon_\Cc:\ \Cc\to S$, satisfying the usual coassociativity and counit conditions:
$$(\Cc\ot_S \Delta_\Cc)\circ \Delta_\Cc=(\Delta_\Cc\ot_S \Cc)\circ \Delta_\Cc~~;~~
(\Cc\ot_S \varepsilon_\Cc)\circ \Delta_\Cc=(\varepsilon_\Cc\ot_S \Cc)\circ \Delta_\Cc=\Cc.$$
For the comultiplication $\Delta_\Cc$, we use the following Sweedler type notation:
$$\Delta_\Cc(c)=c_{(1)}\ot_S c_{(2)}.$$
A right $\Cc$-comodule $M$ is a right $S$-module together with a right $S$-linear map
$\rho:\ M\to M\ot_S \Cc$ such that
$$(M\ot_S\Delta_\Cc)\circ \rho=(\rho\ot_S M)\circ\rho~~;~~
(M\ot_S\varepsilon_\Cc)\circ \rho=S.$$
If $\Cc$ is an $S$-coring, then ${}_S\Hom(\Cc,S)$ is an $S$-ring. This means
that ${}_S\Hom(\Cc,S)$ is a ring, and that we have a ring morphism $j:\ S\to {}_S\Hom(\Cc,S)$. The multiplication
on ${}_S\Hom(\Cc,S)$ is given by the formula
\begin{equation}\eqlabel{left}
(g\# f)(c)=f(c_{(1)}g(c_{(2)})).
\end{equation}
The unit is $\varepsilon_\Cc$, and $j(s)(c)=\varepsilon_\Cc(c)s$, for all $s\in S$ and $c\in\Cc$. In a similar way, $\Hom_S(\Cc,S)$ is an $S$-ring. The multiplication
is now given by the formula
\begin{equation}\eqlabel{right}
(f\# g)(c)=f(g(c_{(1)})c_{(2)}).
\end{equation}

For a detailed discussion of corings and their applications, we refer to \cite{BrzezinskiWisbauer}.

Let $S$ be a commutative $R$-algebra. We have seen in \seref{2.1} that we have
a functor $G:\ \Mm_{S^{\ot 2}}\to {}_S\Mm_S$. An $S$-bimodule $M$ lies in the image
of $G$ if $M^R=M$, that is, $xm=mx$, for all $m\in M$ and $x\in R$.\\
We can view $\Mm_{S^{\ot 2}}$
as a monoidal category with tensor product $\ot_S$ and unit object $S$. A coalgebra
in this category will be called an $S/R$-coring. Thus an $S/R$-coring $\Cc$ is an
$S$-coring, with the additional condition that $\Cc^R=\Cc$.

\begin{example}\exlabel{2.3}
Take an invertible $S$-module $I$. Then $I$ is finitely projective as an $S$-module,
and we have a finite dual basis $\{(e_i,f_i)\in I\times {I^*}~|~i=1,\cdots,n\}$ of $I$. Then
$\sum_i e_i\ot_S f_i=1\in I\ot_S I^*\cong S$. We have an $S/R$-coring
$$\Cc=\Can_R(I;S)=I^*\ot_R I,$$
with structure maps
$$\Delta_\Cc:\ I^*\ot_R I\to I^*\ot_R I\ot_SI^*\ot_R I\cong I^*\ot_RS\ot_R I$$
$$\varepsilon_\Cc:\ I^*\ot_R I\to S$$
given by 
$$\Delta_\Cc(f\ot x)=\sum_i f\ot e_i\ot_S f_i\ot x=f\ot 1\ot x~~;~~
\varepsilon_\Cc(f\ot x)=f(x).$$
We call $\Cc$ an {\sl elementary} coring. If $I=S$, then we obtain
Sweedler's canonical coring, introduced in \cite{Sweedler75}; in general,
$\Can_R(I;S)$ is an example of a comatrix coring, as introduced in \cite{Kaoutit}.
We also compute that
$${}_S\Hom(\Cc,S)={}_S\Hom(I^*\ot_RI,S)\cong {}_R\Hom(I,I)={}_R\End(I).$$
${}_R\End(I)$ is an $R$-algebra (under composition) and an $S$-ring, and we find an isomorphism
of $S$-rings
\begin{equation}\eqlabel{2.3.5}
{}_S\Hom(\Cc,S)\cong {}_R\End(I)^{\rm op}.
\end{equation}
\end{example}

\begin{lemma}\lelabel{2.4}
Let $S$ and $T$ be commutative $R$-algebras. Then we have a strongly
monoidal functor
$$F=-\ot_R T:\ \Mm_{S\ot_R S}\to \Mm_{(S\ot_R T)\ot_T(S\ot_R T)}=\Mm_{S\ot_R S\ot_R T}.$$
Consequently, if $\Cc$ is an $S/R$-coring, then $F(\Cc)=\Cc\ot_RT$ is an $S\ot_R T/T$-coring.
\end{lemma}

\begin{proof}
$F(M)=M\ot_R T$ is an $S\ot_R T$-bimodule, via $(s\ot t)\cdot (m\ot t'')\cdot (s'\ot t')=
sms'\ot tt''t'$. $F$ is strongly monoidal since $F(S)=S\ot_R T$ and
\begin{eqnarray*}
&&\hspace*{-2cm}
F(M\ot_S N)=(M\ot_S N)\ot_R T\\
&\cong & (M\ot_R T)\ot_{S\ot_RT}(N\ot_RT)= F(M)\ot_{S\ot_R T}F(N).
\end{eqnarray*}
\end{proof}

\begin{example}\exlabel{2.5}
Let $I$ be an invertible $S$-module. Then
\begin{eqnarray*}
&&\hspace*{-2cm}F(\Can_R(I;S))=(I^*\ot_RI)\ot_R T\cong (I^*\ot_RT)\ot_{R\ot_R T}
(I\ot_RT)\\
&\cong & (I\ot_R T)^*\ot_{R\ot_R T} (I\ot_RT)\cong \Can_{T}(I\ot_R T;S\ot_RT).
\end{eqnarray*}
\end{example}

\section{Azumaya corings}\selabel{3}

\begin{lemma}\lelabel{3.1}
Let $S$ be a commutative faithfully flat $R$-algebra, and $I\in \dul{\Pic}(S\ot S)$.
Consider an $S$-bimodule map $\Delta:\ I\to I\ot_S I$, and assume that its corresponding
map $\tilde{\Delta}:\ I_2\to I_3\ot_{S^{\ot 3}}I_1$ in $\Mm_{S^{\ot 3}}$ is
an isomorphism. Then we have an isomorphism of $S^{\ot 3}$-modules
\begin{equation}\eqlabel{3.1.1}
\alpha^{-1}=(\tilde{\Delta}\ot_{S^{\ot 3}}I^*_2)\circ \coev_{I_2}:\
S^{\ot 3}\to I_2\ot_{S^{\ot 3}} I^*_2\to I_3\ot_{S^{\ot 3}}I_1\ot_{S^{\ot 3}} I^*_2.
\end{equation}
$\Delta$ is coassociative if and only if 
$(I,\alpha)\in \dul{Z}^1(S/R,\dul{\Pic})$.
\end{lemma}

\begin{proof}
We have the following isomorphisms of $S^{\ot 4}$-modules:
$$\tilde{\Delta}_1:\ I_{21}=I_{13}\to I_{31}\ot_{S^{\ot 4}} I_{11}=
I_{14}\ot_{S^{\ot 4}} I_{12};$$
$$\tilde{\Delta}_2:\ I_{22}=I_{23}\to I_{32}\ot_{S^{\ot 4}} I_{12}=
I_{24}\ot_{S^{\ot 4}} I_{12};$$
$$\tilde{\Delta}_3:\ I_{23}\to I_{33}\ot_{S^{\ot 4}} I_{13}=
I_{34}\ot_{S^{\ot 4}} I_{13};$$
$$\tilde{\Delta}_4:\ I_{24}\to I_{34}\ot_{S^{\ot 4}} I_{14}.$$
$(I,\alpha)\in \dul{Z}^1(S/R,\dul{\Pic})$ if and only if
the composition
$$
\begin{matrix}
I_{23}&\rTo^{I_{23}\ot \coev_{I_{13}}}&
I_{23}\ot_{S^{\ot 4}}I^*_{13}\ot_{S^{\ot 4}}I_{13}\\
&\rTo^{\tilde{\Delta}_3\ot I^*_{13}\ot \tilde{\Delta}_1}&
I_{34}\ot_{S^{\ot 4}}I_{13}\ot_{S^{\ot 4}}I^*_{13}\ot_{S^{\ot 4}}
I_{14}\ot_{S^{\ot 4}}I_{12}\\
&\rTo^{I_{34}\ot\ev_{I_{13}}\ot I_{14}\ot I_{12}}&
I_{34}\ot_{S^{\ot 4}}
I_{14}\ot_{S^{\ot 4}}I_{12}
\end{matrix}
$$
equals the composition
$$\begin{matrix}
I_{23}&\rTo^{\coev_{I_{24}}\ot I_{23}}&
I_{24}\ot_{S^{\ot 4}} I^*_{24}\ot_{S^{\ot 4}}I_{23}\\
&\rTo^{\tilde{\Delta}_4\ot I^*_{24}\ot\tilde{\Delta}_2}&
I_{34}\ot_{S^{\ot 4}} I_{14}\ot_{S^{\ot 4}}I^*_{24}\ot_{S^{\ot 4}}I_{24}
\ot_{S^{\ot 4}}I_{12}\\
&\rTo^{I_{34}\ot I_{14}\ot \ev_{I_{24}}\ot I_{12}}&
I_{34}\ot_{S^{\ot 4}}  I_{14}\ot_{S^{\ot 4}} I_{12}.
\end{matrix}$$
Let $\{(e_i,e_i^*)~|~i=1,\cdots,n\}$ be a finite dual basis of $I$.
For all $c\in I$, we compute
\begin{eqnarray*}
&&\hspace*{-15mm}
\Bigl(\Bigl({I_{34}\ot\ev_{I_{13}}\ot I_{14}\ot I_{12}}\Bigr)
\circ \Bigl(\tilde{\Delta}_3\ot I^*_{13}\ot \tilde{\Delta}_1 \Bigr)
\circ \Bigl(I_{23}\ot \coev_{I_{13}} \Bigr)\Bigr)(c_{23})\\
&=& \Bigl(\Bigl({I_{34}\ot\ev_{I_{13}}\ot I_{14}\ot I_{12}}\Bigr)
\circ \Bigl(\tilde{\Delta}_3\ot I^*_{13}\ot \tilde{\Delta}_1 \Bigr)\Bigr)
\bigl(\sum_i c_{23}\ot e^*_{i13}\ot e_{i13}\bigr)\\
&=& \Bigl({I_{34}\ot\ev_{I_{13}}\ot I_{14}\ot I_{12}}\Bigr)
\Bigl(\sum_i c_{(1)34}\ot c_{(2)13}\ot e^*_{i13}\ot \tilde{\Delta}_1(e_{i13})\Bigr)\\
&=& \sum_i c_{(1)34}\ot\tilde{\Delta}_1((\lan c_{(2)},e_i^*\ran e_i)_{13})
= c_{(1)34}\ot \tilde{\Delta}_1(c_{(2)13})\\
&=& c_{(1)34}\ot c_{(2)(1)14}\ot c_{(2)(2)12},
\end{eqnarray*}
and
\begin{eqnarray*}
&&\hspace*{-15mm}
\Bigl(\Bigl(I_{34}\ot I_{14}\ot \ev_{I_{24}}\ot I_{12} \Bigr)
\circ \Bigl(\tilde{\Delta}_4\ot I^*_{24}\ot\tilde{\Delta}_2\Bigr)
\circ \Bigl(\coev_{I_{24}}\ot I_{23}\Bigr)\Bigr)(c_{23})\\
&=&
\Bigl(\Bigl(I_{34}\ot I_{14}\ot \ev_{I_{24}}\ot I_{12} \Bigr)
\circ \Bigl(\tilde{\Delta}_4\ot I^*_{24}\ot\tilde{\Delta}_2\Bigr)\Bigr)
\bigl(\sum_i e_{i24}\ot e^*_{i24}\ot c_{23}\bigr)\\
&=& \Bigl(I_{34}\ot I_{14}\ot \ev_{I_{24}}\ot I_{12} \Bigr)
\Bigl(\sum_i \tilde{\Delta}_4(e_{i24})\ot e^*_{i24}\ot c_{(1)24}\ot c_{(2)12}\Bigr)\\
&=& \tilde{\Delta}_4(c_{(1)24})\ot c_{(2)12}
= c_{(1)(1)34}\ot c_{(1)(2)14}\ot c_{(2)12}.
\end{eqnarray*}
From \leref{2.1a}, it follows that $(I,\alpha)\in \dul{Z}^1(S/R,\dul{\Pic})$
if and only if the maps in $\Hom_{S^{\ot 4}}(I_{23}, I_{34}\ot_{S^{\ot 4}}  I_{14}\ot_{S^{\ot 4}} I_{12})$ associated to 
$(\Delta\ot_SI)\circ \Delta$ and $(I\ot_S\Delta)\circ \Delta$
in ${}_S\Hom_S(I,I\ot_SI\ot_S I)$ are equal. This is equivalent to the coassociativity
of $\Delta$.
\end{proof}

Observe that the map $\tilde{\Delta}$ can be recovered from $\alpha$ using the
following formula
\begin{equation}\eqlabel{3.1.2}
\tilde{\Delta}=(I_3\ot I_1\ot \ev_{I_2})\circ (\alpha^{-1}\ot I_2).
\end{equation}

\begin{lemma}\lelabel{3.2}
Let $I,\Delta,\tilde{\Delta},\alpha$ be as in \leref{3.1}, and take
$J\in \dul{\Pic}(S)$. Then we have an isomorphism of bimodules with coassociative
comultiplication
$I\cong\Can_R(J;S)$ if and only if $(I,\alpha)
\cong \delta_0(J)$ in $\dul{Z}^1(S/R,\dul{\Pic})$.
\end{lemma}

\begin{proof}
Take $J\in \Pic(S)$. Then $\delta_0(J)=J_1\ot_{S^{\ot 2}} J^*_2=
J^*\ot J$,  and
\begin{eqnarray*}
&&\hspace*{-2cm}
\delta_1\delta_0(J)=\delta_0(J)_1 \ot_{S^{\ot 3}} \delta_0(J)_3 \ot_{S^{\ot 3}} \delta_0(J^*)_2\\
&=&
J_{11}\ot_{S^{\ot 3}}J^*_{21}\ot_{S^{\ot 3}}
J_{13}\ot_{S^{\ot 3}}J^*_{23}\ot_{S^{\ot 3}}
J^*_{12}\ot_{S^{\ot 3}}J_{22}\\
&=&
J_{12}\ot_{S^{\ot 3}}J^*_{13}\ot_{S^{\ot 3}}
J_{13}\ot_{S^{\ot 3}}J^*_{23}\ot_{S^{\ot 3}}
J^*_{12}\ot_{S^{\ot 3}}J_{23}.
\end{eqnarray*}
The map $\lambda_J$ is obtained by applying the evaluation map
on tensor factors 1 and 5, 2 and 3, 4 and 6. Let $\{(e_i,e_i^*)~|~i=1,\cdots, n\}$
be a finite dual basis of $J$ as an $S$-module. Then
$$\lambda_J^{-1}(1\ot 1\ot 1)=
\sum_{i,j,k} e_{i12}\ot e^*_{j13}\ot e_{j13}\ot e^*_{k23}\ot e^*_{i12}\ot e_{k23}.$$
Take $x^*\ot 1\ot x=x_{12}\ot_{S^{\ot 3}} x_{23}^*\in (J^*\ot J)_2=J^*\ot S\ot J=J_{12}\ot_{S^{\ot 3}}J^*_{23}$. We then compute, using \equref{3.1.2},
\begin{eqnarray*}
&&\hspace*{-12mm}
\tilde{\Delta}(x^*\ot 1\ot x)=
(\delta_0(J)_1 \ot_{S^{\ot 3}} \delta_0(J)_3 \ot_{S^{\ot 3}} \ev_{\delta_0(J)_2})
(\lambda_J^{-1}\ot_{S^{\ot 3}} \delta_0(J)_2)(x_{12}\ot x_{23}^*)\\
&=&
\sum_{i,j,k} e_{i12}\ot e^*_{j13}\ot e_{j13}\ot e^*_{k23}\ot 
\lan e^*_{i12}, x_{12}\ran \ot \lan x_{23}^*,e_{k23}\ran\\
&=& \sum_{j} x_{12}\ot e_{j13}\ot e^*_{j13}\ot x^*_{23}\\
&=& \sum_j x^*\ot (e^*_j\ot_S e_j)\ot x\in J^*\ot (J^*\ot_S J)\ot J.
\end{eqnarray*}
Consequently
$$\Delta(x^*\ot x)=\sum_j x^*\ot \lan e^*_j, e_j\ran \ot x=x^*\ot 1\ot x$$
is the comultiplication on $\Can_R(J;S)$.
In a similar way, starting from the comultiplication $\Delta$ on
$\Can_R(J;S)$, we find that the map $\alpha$ defined in \equref{3.1.1}
is precisely $\lambda_J$.
\end{proof}

\begin{theorem}\thlabel{3.3}
Let $\Cc$ be a faithfully projective $S\ot S$-module, and $\Delta:\
\Cc\to \Cc\ot_S\Cc$ an $S$-bimodule map. We consider the corresponding map
$\tilde{\Delta}:\ \Cc_2\to \Cc_3\ot_{S^{\ot 3}}\Cc_1$ in $\Mm_{S^{\ot 3}}$
(cf. \leref{2.1}). Then the following assertions are equivalent.
\begin{enumerate}
\item $\Delta$ is coassociative and $\tilde{\Delta}$ is an isomorphism in
$\Mm_{S^{\ot 3}}$;
\item $\Cc\in \dul{\Pic}(S^{\ot 2})$ and $(\Cc,\alpha)\in \dul{Z}^1(S/R,
\dul{\Pic})$, with $\alpha$ defined by \equref{3.1.1};
\item $\Cc\in \dul{\Pic}(S^{\ot 2})$ and $\Cc\ot S$ is isomorphic to
$\Can_{R\ot S}(\Cc;S\ot S)$ as bimodules with coassociative comultiplication;
\item there exists a faithfully flat commutative $R$-algebra $T$
such that $(\Cc\ot_R T,\Delta\ot_R T)$ is isomorphic to 
$\Can_T(I;S\ot T)$, for some $I\in \dul{\Pic}(S\ot T)$, as a bimodule with a coassociative
comultiplication;
\item $(\Cc,\Delta)$ is a coring and $\tilde{\Delta}$ is an isomorphism in
$\Mm_{S^{\ot 3}}$.
\end{enumerate}
\end{theorem}

\begin{proof}
$\ul{ 1)\Rightarrow 2)}$. From the fact that $\tilde{\Delta}$ is
an isomorphism, it follows that $\Cc_2\cong \Cc_3\ot_{S^{\ot 3}}\Cc_1$.
Applying the functor $|-|:\ \Mm_{S^{\ot 3}}\to \Mm_S$, we find that
$|\Cc|\cong |\Cc|\ot_S |\Cc|$. $\Cc$ is a faithfully projective $S\ot S$-module,
so $|\Cc|$ is a faithfully projective $S$-module. Its rank is an idempotent,
so it is equal to one, and $|\Cc|\in \dul{\Pic}(S)$.\\
Now switch the second and third tensor factor in $\Cc_2\cong \Cc_3\ot_{S^{\ot 3}}\Cc_1$,
and then apply $|-|$ to the first and second factor. We find that
$|\Cc|\ot S\cong \Cc\ot_{S^{\ot 2}} \tau(\Cc)$, with $\tau(\Cc)$ equal to $\Cc$
as an $R$-module, with newly defined $S\ot S$-action $c\triangleleft(s\ot t)=
c(t\ot s)$. Now $|\Cc|\ot S\in \dul{\Pic}(S\ot S)$, and it follows that
$\Cc\in \dul{\Pic}(S\ot S)$. It follows now from \leref{3.1} that
$(\Cc,\alpha)\in \dul{Z}^1(S/R,\dul{\Pic})$.\\
$\ul{ 2)\Rightarrow 3)}$. It follows from \leref{1.2a.2} that
$(\Cc\ot S,\alpha\ot S)\cong \delta_0(\Cc)$ in $\dul{Z}^1(S\ot S/R\ot S,\dul{\Pic})$.
From \leref{3.2}, it follows that $\Cc\ot S\cong
\Can_{R\ot S}(\Cc;S\ot S)$ as bimodules with coassociative comultiplication.\\
$\ul{ 3)\Rightarrow 4)}$ is obvious.\\
$\ul{ 4)\Rightarrow 1)}$. After faithfully flat base extension, $\Delta$ becomes
coassociative, and $\tilde{\Delta}$ becomes an isomorphism. Hence $\Delta$ is
coassociative and $\tilde{\Delta}$ is an isomorphism.\\
$\ul{ 1)\Rightarrow 5)}$. We have an isomorphism of $S^{\ot 3}$-modules
$\alpha:\ \Cc^*_2\ot_{S^{\ot 3}}\Cc_1\ot_{S^{\ot 3}} \Cc_3\to S^{\ot 3}$. Applying the functor
$|-|$, we find an isomorphism of $S$-modules
$|\alpha |:\ |\Cc|\to S$. Now we consider the composition $\varepsilon=
|\alpha|\circ \pi:\ \Cc\to S$. In the situation where $\Cc=\Can_R(I;S)$, $\varepsilon$
is the counit of $\Cc$. By 4), $\varepsilon$ has the counit property after
a base extension. Hence $\varepsilon$ has itself the counit property.
So $(\Cc,\Delta,\varepsilon)$ is a coring.\\
$\ul{5)\Rightarrow 1)}$ is obvious.
\end{proof}

If $(\Cc,\Delta,\varepsilon)$ satisfies the equivalent conditions of
\thref{3.3}, then we call $\Cc$ an Azumaya $S/R$-coring. The connection to
Azumaya algebras is discussed in the following Proposition.

\begin{proposition}\prlabel{3.3a}
Let $S$ be a faithfully projective commutative $R$-algebra, and $\Cc$
an Azumaya $S/R$-coring. Then ${}_S\Hom(\Cc,S)$ and $\Hom_S(\Cc,S)$ are Azumaya $R$-algebras
split by $S$.
\end{proposition}

\begin{proof}
Using \thref{3.3} and \equref{2.3.5}, we find the following isomorphisms
of $S$-algebras:
\begin{eqnarray*}
&&\hspace*{-2cm}
{}_S\Hom(\Cc,S)\ot S={}_S\Hom(\Cc,S)\ot {}_S\Hom(S,S)\cong 
{}_{S\ot S}\Hom(\Cc\ot S,S\ot S)\\
&\cong & {}_{S\ot S}\Hom(\Can_{R\ot S}(\Cc;S\ot S),S\ot S)\cong
{}_{R\ot S}\End(\Cc)^{\rm op}.
\end{eqnarray*}
\end{proof}

\begin{theorem}\thlabel{3.4}
Let $(\Cc,\Delta)$ and $(\Cc',\Delta')$ be Azumaya $S/R$-corings, and 
consider the corresponding $(\Cc,\alpha),(\Cc',\alpha')\in
\dul{Z}^1(S/R,\dul{\Pic})$. Let $f:\ \Cc\to \Cc'$ be an isomorphism
in $\dul{\Pic}(S\ot S)$. Then $f$ is an isomorphism of corings if and only
if $f$ defines an isomorphism in $\dul{Z}^1(S/R,\dul{\Pic})$.
\end{theorem}

\begin{proof}
$f$ is an isomorphism of corings if and only if the following diagram
commutes:
$$\begin{diagram}
\Cc&\rTo^{\Delta}&\Cc\ot_S\Cc\\
\dTo_{f}&&\dTo_{f\ot_S f}\\
\Cc'&\rTo^{\Delta'}&\Cc'\ot_S\Cc'.
\end{diagram}$$
This is equivalent to commutativity of the diagram
$$\begin{diagram}
\Cc_2&\rTo^{\tilde{\Delta}}&\Cc_3\ot_{S^{\ot 3}}\Cc_1\\
\dTo_{f_2}&&\dTo_{f_3\ot_{S^{\ot 3}} f_1}\\
\Cc'_2&\rTo^{\tilde{\Delta}'}&\Cc'_3\ot_{S^{\ot 3}}\Cc'_1.
\end{diagram}$$
This is equivalent to commutativity of the right square in the next diagram
$$\begin{diagram}
S^{\ot 3}&\rTo^{\coev_{\Cc_2}}&
\Cc^*_2\ot_{S^{\ot 3}}\Cc_2&\rTo^{\Cc^*_2\ot\tilde{\Delta}}&\Cc^*_2\ot_{S^{\ot 3}}\Cc_3\ot_{S^{\ot 3}}\Cc_1\\
\dTo_{=}&&
\dTo_{(f^*_2)^{-1}\ot_{S^{\ot 3}} f_2}&&\dTo_{(f^*_2)^{-1}\ot_{S^{\ot 3}} f_3\ot_{S^{\ot 3}} f_1}\\
S^{\ot 3}&\rTo^{\coev_{\Cc'_2}}&
{\Cc'}^*_2\ot_{S^{\ot 3}}\Cc'_2&\rTo^{{\Cc'}^*_2\ot\tilde{\Delta}'}&{\Cc'}^*_2\ot_{S^{\ot 3}}\Cc'_3
\ot_{S^{\ot 3}}\Cc'_1.
\end{diagram}$$
The left square is automatically commutative. Commutativity of the full diagram
is equivalent to $\alpha'\circ \delta_1(f)=\alpha$, as needed.
\end{proof}

Let ${\Az}^c(S/R)$ be the category of Azumaya $S/R$-corings and isomorphisms of corings.

\begin{proposition}\prlabel{3.5}
$(\Az^c(S/R),\ot_{S^{\ot 2}},\Can_R(S;S))$ is a monoidal category.
\end{proposition}

\begin{proof}
Take two Azumaya $S/R$-corings $(\Cc,\Delta)$ and $(\Cc',\Delta')$, and let
$\tilde{D}$ be the following composition
\begin{eqnarray*}
(\Cc\ot_{S^{\ot 2}}\Cc')_2 = \Cc_2\ot_{S^{\ot 3}}\Cc'_2
&\rTo^{\tilde{\Delta}\ot \tilde{\Delta}'}&
\Cc_3\ot_{S^{\ot 3}}\Cc_1 \ot_{S^{\ot 3}} \Cc'_3\ot_{S^{\ot 3}}\Cc'_1\\
&\rTo^{\Cc_3\ot \tau\ot \Cc'_1}&
\Cc_3\ot_{S^{\ot 3}}\Cc'_3 \ot_{S^{\ot 3}} \Cc_1\ot_{S^{\ot 3}}\Cc'_1\\
&=& (\Cc\ot_{S^{\ot 2}}\Cc')_3 \ot_{S^{\ot 3}}(\Cc\ot_{S^{\ot 2}}\Cc')_1.
\end{eqnarray*}
The comultiplication on $\Cc\ot_{S^{\ot 2}}\Cc'$ is the corresponding map
$$D:\ \Cc\ot_{S^{\ot 2}}\Cc' \to (\Cc\ot_{S^{\ot 2}}\Cc')\ot_S (\Cc\ot_{S^{\ot 2}}\Cc').$$
Observe that the $S$-bimodule structure on $\Cc\ot_{S^{\ot 2}}\Cc'$ is
given by the formulas
$$s(c\ot c')=sc\ot c'=c\ot sc'~~;~~
(c\ot c')t=c\ot c' t= ct\ot c'.$$
We have that
$$\tilde{D}(c\ot c')_2= (c_{(1)}\ot_{S^{\ot 2}} c'_{(1)})_{3}
\ot_{S^{(3)}} (c_{(2)}\ot_{S^{\ot 2}} c'_{(2)})_{1},$$
hence
$$D(c\ot c')=(c_{(1)}\ot_{S^{\ot 2}} c'_{(1)})
\ot_S (c_{(2)}\ot_{S^{\ot 2}} c'_{(2)}).$$
It is then easy to see that $D$ is coassociative, and that
$$\Cc\ot_{S^{\ot 2}} \Can_R(S;S)\cong \Cc\cong \Can_R(S;S)\ot_{S^{\ot 2}}\Cc.$$
\end{proof}

\begin{corollary}\colabel{3.6}
We have a monoidal isomorphism of categories
$$H:\ \Az^c(S/R)\to \dul{Z}^1(S/R,\dul{\Pic}).$$
\end{corollary}

Consider the subgroup $\Can^c(S/R)$ of $K_0\Az^c(S/R)$ consisting of
isomorphism classes represented by an elementary coring $\Can_R(I;S)$ for some $I\in \dul{\Pic}(S)$.
The quotient
$$\Br^c(S/R)=K_0\Az^c(S/R)/\Can^c(S/R)$$
is called the {\sl relative Brauer group} of Azumaya $S/R$-corings.

\begin{corollary}\colabel{3.7}
We have an isomorphism of abelian groups
$$\Br^c(S/R)\cong H^1(S/R,\dul{\Pic}).$$
Consequently, we have an exact sequence
\begin{eqnarray}\eqlabel{3.7.1}
0&\longrightarrow& H^1(S/R,\units)\longrightarrow \Pic(R)
\longrightarrow H^0(S/R,\Pic)\\
&\longrightarrow& H^2(S/R,\units)\rTo^{\alpha} \Br^c(S/R)
\longrightarrow H^1(S/R,\Pic)\nonumber\\
&\longrightarrow&  H^{3}(S/R,\units)\nonumber.
\end{eqnarray}
\end{corollary}

Let $f:\ S\to T$ be a morphism of faithfully flat commutative $R$-algebras.
Then we have a functor $\tilde{f}:\ \Az^c(S/R)\to \Az^c(T/R)$ such that the
following diagram commutes
$$\begin{diagram}
\Az^c(S/R)&\rTo^{H}& \dul{Z}^1(S/R,\dul{\Pic})\\
\dTo^{\tilde{f}}&&\dTo_{f_*}\\
\Az^c(T/R)&\rTo^{H}& \dul{Z}^1(T/R,\dul{\Pic}).
\end{diagram}$$
$\tilde{f}(\Cc)=\Cc\ot_{S^{\ot 2}}\Can_R(T;T)$, with comultiplication
$\Delta_\Cc\ot_{S^{\ot 2}}\Delta$, where $\Delta$ is the comultiplication on
the canonical coring $\Can_R(T;T)$. This induces a commutative diagram
$$\begin{diagram}
\Br^c(S/R)&\rTo^{\cong}& H^1(S/R,\dul{\Pic})\\
\dTo^{\tilde{f}}&&\dTo_{f_*}\\
\Br^c(T/R)&\rTo^{\cong}& H^1(T/R,\dul{\Pic}).
\end{diagram}$$
Otherwise stated, the isomorphisms in
\coref{3.7} define an isomorphism of functors
$$\Br^c(\bullet/R)\cong H^1(\bullet/R,\dul{\Pic}):\ \Rr\to\dul{Ab},$$
and
\begin{equation}\eqlabel{3.7.2}
\colim \Br^c(\bullet/R)\cong \check H^1(R_{\rm fl},\dul{\Pic})
\cong H^2(R_{\rm fl},\units).
\end{equation}

Let us describe the map $\alpha:\ H^2(S/R,\units)\to \Br^c(S/R)$.
Let $u\in Z^2(S/R,\units)$ be a cocycle, and consider the coring
$$\Cc=\Can_R(S;S)_u,$$
which is equal to $S\ot S$ as an $S$-bimodule, with comultiplication
\begin{equation}\eqlabel{3.7.3}
\Delta_u:\ S\ot S\to S\ot S\ot_S S\ot S\cong S\ot S\ot S,~~
\Delta_u(s\ot t)=u^1s\ot u^2\ot u^3t.
\end{equation}
The coassociativity follows immediately from the cocycle condition;
the counit $\varepsilon$ is given by the formula (see \leref{am3})
\begin{equation}\eqlabel{3.7.4}
\varepsilon(s\ot t)=|u|^{-1}st.
\end{equation}
The counit property follows from \leref{am3}. If $u$ is normalized, then
the counit coincides with the counit in $\Can_R(S;S)$.\\
Let us compute the right dual $\Hom_S(\Cc,S)$. As an $R$-module,
$\Hom_S(\Cc,S)=\Hom_S(S\ot S,S)\cong \End_R(S)$. We transport the multiplication on
$\Hom_S(\Cc,S)$ to $\End_R(S)$ as follows: take $\varphi,\psi\in \End_R(S)$, and
define $f,g\in \Hom_S(S\ot S,S)$ by
$$f(s\ot t)=\varphi(s)t~~;~~g(s\ot t)=\psi(s)t.$$
Then we find, using \equref{right},
$$(\varphi * \psi)(s)=(f\# g)(s\ot 1)=f(\psi(su^1)u^2\ot u^3)=
\varphi(\psi(su^1)u^2)u^3,$$
or
\begin{equation}\eqlabel{3.8.1}
\varphi*\psi=u^3\varphi u^2\psi u^1.
\end{equation}
In a similar way, we find that ${}_S\Hom(\Cc,S)\cong \End_R(S)$, with twisted
multiplication
\begin{equation}\eqlabel{3.8.2}
\varphi*\psi=u^1\psi u^2\varphi u^3.
\end{equation}

If $S$ is faithfully projective as an $R$-module, then it is
well-known that there exists a morphism
$$\alpha:\ H^2(S/R, \units)\to \Br(S/R).$$
More precisely, we can associate an Azumaya algebra $A(u)$ to any cocycle
$u\in Z^2(S/R, \units)$. The construction of $A(u)$ was given first
in \cite[Theorem 2]{RZ}. It is explained in \cite[V.2]{KO1} and
\cite[7.5]{K2} using descent theory. Let us summarize the construction of
$A(u)$, following \cite{KO1}. Take a cocycle $u=u^1\ot u^2\ot u^3=
U^1\ot U^2\ot U^3$ with inverse $u^{-1}= v^1\ot v^2\ot v^3$, and
consider the map
$$\Phi:\ S\ot S\ot \End_R(S)\to S\ot \End_R(S)\ot S,~~
\Phi(s\ot t\ot \varphi)=su^1v^1\ot u^3\varphi v^3 \ot tu^2v^2.$$
Then
$$A(u)=\{x\in S\ot \End_R(S)~|~x\ot 1=\Phi(1\ot x)\}.$$
It will be convenient to use the canonical identification $\End_R(S)\cong S^*\ot S$.
Then $x=\sum_i s_i\ot t^*_i\ot t_i\in S\ot S^*\ot S$ lies in $A(u)$ if and only
if
$$\sum_i s_i\ot t_i^*\ot t_i\ot 1=
\sum_i u^1v^1\ot t_i^*v^3\ot u^3t_i\ot u^2v^2s_i,$$
or
$$\sum_i s_i\ot 1\ot t_i^*\ot t_i=
\sum_i u^1v^1\ot u^2v^2s_i\ot t_i^*v^3\ot u^3t_i,$$
or
\begin{equation}\eqlabel{3.8.3}
x_2=x_1u_3u_4^{-1}~~{\rm or}~~x_2u_4=x_1u_3.
\end{equation}

Let $\End_R(S)_u$ be equal to $\End_R(S)$, with twisted multiplication given
by \equref{3.8.1}. We know from \prref{3.3a} that $\End_R(S)_u$ is an
Azumaya algebra split by $S$.

\begin{theorem}\thlabel{3.8}
Let $S$ be a faithfully projective commutative $R$-algebra, and $u\in Z^2(S/R,\units)$.
Then we have an isomorphism of $R$-algebras $\gamma:\ \End_R(S)_u\to A(u)$.
\end{theorem}

\begin{proof}
We define $\gamma$ by the following formula:
$$\gamma(\varphi)=u^1\ot u^3\varphi u^2,$$
or
$$\gamma(t^*\ot t)=u^1 \ot t^*u^2\ot u^3t.$$
We have to show that $x=\gamma(t^*\ot t)$ satisfies \equref{3.8.3}. Indeed,
$$x_2u_4=(1\ot 1\ot t^*\ot t)u_2u_4=(1\ot 1\ot t^*\ot t)u_1u_3=x_1u_3.$$
Let us next show that $\gamma$ is multiplicative. We want to show that
$$\gamma(\psi)\circ \gamma(\varphi)=\gamma(\psi *\varphi)$$
or
$$u^1U^1\ot u^3\psi u^2 U^3\varphi U^2=
U^1\ot U^3u^3\psi u^2\varphi u^1 U^2.$$
It suffices that
$$u^1U^1\ot u^3\ot u^2 U^3\ot U^2=
U^1\ot U^3u^3\ot u^2\ot u^1 U^2,$$
or
$$u^1U^1\ot U^2\ot u^2 U^3\ot u^3=
U^1\ot u^1 U^2\ot u^2\ot U^3u^3.$$
This is precisely the cocycle condition $u_2u_4=u_1u_3$.\\
The inverse of $\gamma$ is given by
$$\gamma^{-1}(\sum_i s_i\ot t_i^*\ot t_i)=\sum_i t_i^*v^2\ot v^1v^3s_it_i,$$
for all $x=\sum_i s_i\ot t_i^*\ot t_i\in A(u)$. We compute that
$$\gamma(\gamma^{-1}(x))=
\gamma(\sum_i t_i^*v^2\ot v^1v^3s_it_i)
= u^1\ot t_i^*v^2u^2\ot u^3v^1v^3s_it_i.$$
It follows from \equref{3.8.3} that
$$x_2=x_1u_3u_4^{-1}=x_1u_2u_1^{-1}=
u^1\ot s_iv^1\ot t_i^*u^2v^2\ot t_iu^3v^3.$$
Multiplying the second and the fourth tensor factor, we obtain that
$$\gamma(\gamma^{-1}(x))=u^1\ot t_i^*v^2u^2\ot u^3v^1v^3s_it_i=x.$$
Finally
$$
\gamma^{-1}(\gamma(t^*\ot t))=\gamma^{-1}(u^1\ot t^*u^2\ot u^3t)=
 t^*u^2v^2\ot v^1v^3u^1u^3t =t^*\ot t.$$
\end{proof}

\section{A Normal Basis Theorem}\selabel{normal}
Let $S$ be a faithfully flat commutative $R$-algebra.
We say that an $S\ot S$-module with coassociative comultiplication
has normal basis if it is isomorphic to $S\ot S$ as an $S$-bimodule.
Examples are the Azumaya $S/R$-corings $\Can_R(S;S)_u$, with
$u\in Z^2(S/R,\units)$, as considered above.
The category of $S/R$-corings (resp. $S\ot S$-modules with coassociative
comultiplication)
with normal basis will be denoted by $\dul{F}(S/R)$ (resp. $\dul{F}'(S/R)$).
$(\dul{F}(S/R),\ot_{S^{\ot 2}},\Can_R(S;S))$ and $(\dul{F}'(S/R),\ot_{S^{\ot 2}},$ $\Can_R(S;S))$
are monoidal categories, and the
sets of isomorphism classes $F(S/R)$ and $F'(S/R)$ are monoids. Let
$FAz(S/R)$ be the subgroup of $F(S/R)$ consisting of isomorphism classes of
$S/R$-Azumaya corings with normal basis.
We have inclusions
$$FAz(S/R)\subset F(S/R)\subset F'(S/R).$$
We will
give a cohomological description of these monoids.\\
Take $u=u^1\ot u^2\ot u^3\in S^{\ot 3}$. As usual, summation is implicitly
understood. We do not assume that $u$ is invertible. We call $u$ a $2$-cosickle
if $u_1u_3=u_2u_3$. If, in addition, $u^1u^2\ot u^3$ and
$u^1\ot u^2u^3$ are invertible in $S^{\ot 2}$, then we call $u$ an
almost invertible $2$-cosickle. This implies in particular that
$|u|=u^1u^2u^3$ is invertible in $S$.
Almost invertible $2$-cosickles have been introduced and studied in \cite{Haile}.
Let ${S'}^2(S/R)$ be the set of $2$-cosickles and $S^2(S/R)$ the set of almost
invertible $2$-cosickles. ${S}^2(S/R)$ and ${S'}^2(S/R)$ are 
multiplicative monoids, and we have the following inclusions of monoids:
$$B^2(S/R,\units) \subset Z^2(S/R,\units)\subset S^2(S/R)\subset {S'}^2(S/R)\subset
S^{\ot 3}.$$
We consider the quotient monoids
$$M^2(S/R)=S^2(S/R)/B^2(S/R,\units)~~;~~{M'}^2(S/R)={S'}^2(S/R)/B^2(S/R,\units).$$
$M^2(S/R)$ is called
the second (Hebrew) Amitsur cohomology monoid; the subgroup consisting of
invertible classes is the usual (French) Amitsur cohomology group $H^2(S/R,\units)$
(the Hebrew-French dictionary is explained in detail in \cite{Haile}). We have the
following inclusions:
$$H^2(S/R,\units)\subset M^2(S/R)\subset {M'}^2(S/R).$$

\begin{theorem}\thlabel{4a.1}
Let $S$ be a commutative faithfully flat $R$-algebra. An $S\ot S$-module with coassociative
comultiplication and normal basis is an Azumaya $S/R$-coring 
if and only if it represents an invertible element of $F'(S/R)$.
Furthermore
$$F'(S/R)\cong {M'}^2(S/R),~~F(S/R)\cong M^2(S/R)
~~{\rm and}~~FAz(S/R)\cong H^2(S/R,\units).$$
\end{theorem}

\begin{proof}
We define a map $\alpha':\  {S'}^2(S/R)\to {F'}(S/R)$ as follows:
$\alpha'(u)=\Can_R(S;S)_u$, with comultiplication given by \equref{3.7.3}.
It is easy to see that $\alpha'$ is a map of monoids.
$\alpha'$ is surjective: let $\Cc=S^{\ot 2}$ with a coassociative comultiplication
$\Delta_\Cc$, and take
$$u=u^1\ot u^2\ot u^3=\Delta_\Cc(1\ot 1)\in S\ot S\ot_S S\ot S\cong S^{\ot 3}.$$
From the coassociativity of $\Delta_\Cc$, it follows that $u_1u_3=u_2u_4$,
so $u\in {S'}^2(S/R)$, and $\alpha'(u)=\Cc$.\\
Take $u\in \Ker\alpha'$. We then have a comultiplication preserving
$S$-bimodule isomorphism
$\varphi:\ \Can_R(S;S)\to \Can_R(S;S)_u$. Put $\varphi(1\ot 1)=v\in S^{\ot 2}$.
From the fact that $\varphi$ is an automorphism of $S^{\ot 2}$ as an
$S$-bimodule, it follows that $v^{-1}=\varphi^{-1}(1\ot 1)$.
$\varphi$ preserves comultiplication, so it follows that
\begin{eqnarray*}
&&\hspace*{-1cm}v_1v_3=(\varphi\ot_S\varphi)(\Delta_1(1\ot 1))\\
&=& \Delta_u(\varphi(1\ot 1))=\Delta_u(v)=v^1u^1\ot u^2\ot u^3v^2
= v_2u,
\end{eqnarray*}
hence $u=\delta_1(v)\in B^2(S/R)$. It follows that 
$F'(S/R)\cong {M'}^2(S/R)$ as monoids.\\
If $u\in  {S}^2(S/R)$, then $\alpha'(u)=\Can_R(S;S)_u$ has counit given by
\equref{3.7.4}. Conversely, let $\Cc\in F(S/R)$, and take $u={\alpha'}^{-1}(\Cc)$.
Let $v=\varepsilon_\Cc(1\ot 1)$. Using the counit property and the fact that
$\varepsilon_\Cc$ is a bimodule map, we then compute that
$$1\ot 1=\varepsilon_\Cc(u^1\ot u^2)\ot u^3=u^1vu^2\ot u^3;$$
$$1\ot 1=u^1\ot \varepsilon_\Cc(u^2\ot u^3)=u^1\ot u^2vu^3.$$
It follows that $u^1u^2\ot u^3$ and
$u^1\ot u^2u^3$ are invertible, and that $v=|u|^{-1}$. Hence $u\in S^2(S/R)$,
and it follows that $\alpha'$ restricts to an epimorphism of monoids
$$\alpha:\ {S}^2(S/R)\to {F}(S/R).$$
It is clear that $\Ker\alpha=B^2(S/R,\units)$, and it follows that
$M^2(S/R)\cong F(S/R)$.\\
If $u\in Z^2(S/R,\units)$, then $\alpha'(u)=\Can_R(S;S)_u$ is an Azumaya
$S/R$-coring. Conversely, let $\Cc$ be an Azumaya $S/R$-coring
with normal basis, and
$u={\alpha'}^{-1}(\Cc)$.  Then $[u]$ is invertible
in $M^2(S/R)$, so there exists $v\in S^2(S/R)$ such that $uv\in B^2(S/R)$.
Since every element in $B^2(S/R)$ is invertible in $S^{\ot 3}$, it follows that
$u\in \units(S^{\ot 3})$, and $u\in Z^2(S/R,\units)$. So $\alpha$ restricts
to an epimorphism $\alpha'':\ Z^2(S/R,\units)\to FAz(S/R)$. Clearly
$\Ker\alpha''=B^2(S/R,\units)$, hence $B^2(S/R,\units)\cong FAz(S/R)$.
\end{proof}

\section{The Brauer group}\selabel{5}
An Azumaya coring over $R$ is a pair $(S,\Cc)$, where $S$ is a faithfully
flat finitely presented commutative $R$-algebra, and $\Cc$ is an Azumaya $S/R$-coring.
A morphism between two Azumaya corings $(S,\Cc)$ and $(T,\Dd)$ over $R$
is a pair $(f,\varphi)$, with $f:\ S\to T$ an algebra isomorphism, and
$\varphi:\ \Cc\to\Dd$ an $R$-module isomorphism preserving the bimodule
structure and the comultiplication, that is
$$\varphi(scs')=f(s)\varphi(c)f(s')~~{\rm and}~~
\Delta_\Dd(\varphi(c))=\varphi(c_{(1)})\ot_T\varphi(c_{(2)}),$$
for all $s,s'\in S$ and $c\in \Cc$. The counit is then preserved
automatically. Let $\dul{\rm Az}^c(R)$ be the category of Azumaya corings over $R$.

\begin{lemma}\lelabel{5.1}
Suppose that $S$ and $T$ are commutative $R$-algebras. If $M\in \Mm_{S\ot_R S}$
and $N\in \Mm_{T\ot_R T}$, then $M\ot_R N\in \Mm_{(S\ot_RT)\ot_R(S\ot_RT)}$.\\
If $\Cc$ is an (Azumaya) $S/R$-coring, and $\Dd$ is an (Azumaya) $T/R$-coring, then
$\Cc\ot_R\Dd$ is an (Azumaya) $S\ot_RT/R$-coring.
\end{lemma}

\begin{proof}
The proof of the first two assertions is easy; the structure maps are the
obvious ones. Let us show that $\Cc\ot_R\Dd$ is an Azumaya $S\ot_RT/R$-coring.
\begin{eqnarray*}
&&\hspace*{-15mm}
\Cc\ot_R\Dd\ot_RS\ot_RT\cong \Cc\ot_RS\ot_R\Dd\ot_RT\\
&\cong & \Can_S(I;S\ot_R S)\ot \Can_T(J;T\ot_R T)=(I^*\ot_SI)\ot_R(J^*\ot_TJ)\\
&\cong & (I^*\ot_RJ^*)\ot_{S\ot T}(I\ot_RJ)=\Can_{S\ot T}(I\ot_RJ;
S\ot_RT\ot_RS\ot_RT).
\end{eqnarray*}
\end{proof}

Let $(\Cc,\Delta)$ be an Azumaya $S/R$-coring, and consider the corresponding
$(\Cc,\alpha)\in \dul{Z}^1(S/R,\dul{\Pic})$. Its inverse in
$Z^1(S/R,\dul{\Pic})$ is represented by $(\Cc^*,(\alpha^*)^{-1})$. The corresponding
coring will be denoted by $(\Cc^*,\ol{\Delta})$.

\begin{proposition}\prlabel{5.2}
Let $\Cc$ be an Azumaya $S/R$-coring. Then
$\Cc\ot \Cc^*$ is an elementary coring.
\end{proposition}

\begin{proof}
Consider $H(\Cc)=(\Cc,\alpha)\in \dul{Z}^1(S/R,\dul{\Pic})$, and the maps
$\eta_1,\eta_2:\ S\to S\ot S$. It follows from \prref{1.2a.3} that
$$[\eta_{1*}(\Cc,\alpha)]=[(\Cc\ot S^{\ot 2},\alpha\ot S^{\ot 3})]
= [\eta_{2*}(\Cc,\alpha)]=[(S^{\ot 2}\ot \Cc,S^{\ot 3}\ot \alpha)]$$
in $H^1(S\ot S/R,\dul{\Pic})$. Consequently
$$[H^{-1}(\eta_{1*}(\Cc,\alpha))]=[\Cc\ot \Can_R(S;S)]=
[H^{-1}(\eta_{2*}(\Cc,\alpha))]=[\Can_R(S;S)\ot\Cc]$$
in $\Br^c(S\ot S/R)$. The inverse of $[\Can_R(S;S)\ot\Cc]$ in $\Br^c(S\ot S/R)$
is represented by $\Can_R(S;S)\ot\Cc^*$. It follows that
$$(\Cc\ot \Can_R(S;S))\ot_{S^{\ot 4}} (\Can_R(S;S)\ot\Cc^*)\cong \Cc\ot \Cc^*$$
is an elementary coring.
\end{proof}

Let $(S,\Cc)$ and $(T,\Dd)$ be Azumaya corings over $R$.
We say that $\Cc$ and $\Dd$ are Brauer equivalent (notation: $\Cc\sim \Dd$)
if there exist elementary corings $\Ee_1$ and $\Ee_2$ over $R$ such that
$\Cc\ot \Ee_1\cong \Dd\ot \Ee_2$ as Azumaya corings over $R$. Since the
tensor product of two elementary corings is elementary, it is easy to show
that 
$\sim$ is an equivalence relation.
Let $\Br^{\rm c}_{\rm fl}(R)$ be the set of equivalence classes of
isomorphism classes of Azumaya corings over $R$.

\begin{proposition}\prlabel{5.4}
$\Br^{\rm c}_{\rm fl}(R)$ is an abelian group under the operation induced by
the tensor product $\ot_R$, with unit element $[R]$.
\end{proposition}

\begin{proof}
It follows from \prref{5.2} that the inverse of $[(\Cc,\Delta)]$ is
$[(\Cc^*,\ol{\Delta})]$.
\end{proof}

\begin{lemma}\lelabel{5.5}
Let $\Cc,\Ee$ be Azumaya $S/R$-corings, and assume that $\Ee=\Can_R(J;S)$ is elementary.
Then the Azumaya corings $\Cc\ot_{S^{\ot 2}}\Ee$ and $\Cc$ are Brauer
equivalent.
\end{lemma}

\begin{proof}
Let $H(\Cc)=(\Cc,\alpha)$. We know that $H(\Ee)=(J^*\ot J,\lambda_J)$, and
$$[(\Cc\ot_{S^{\ot 2}}\Ee,\alpha\ot_{S^{\ot 3}}\lambda_J)]=
[(\Cc,\alpha)]$$
in $H^1(S/R,\dul{\Pic})$. From \prref{1.2a.3}, it follows that
\begin{eqnarray*}
&&\hspace*{-2cm}
[\eta_{1*}(\Cc,\alpha)]=[(S^{\ot 2}\ot\Cc,S^{\ot 3}\ot\alpha)]=
[\eta_{2*}(\Cc\ot_{S^{\ot 2}}\Ee,\alpha\ot_{S^{\ot 3}}\lambda_J)]\\
&=& 
[(\Cc\ot_{S^{\ot 2}}\Ee)\ot S^{\ot 2}, (\alpha\ot_{S^{\ot 3}}\lambda_J)\ot S^{\ot 3})]
\end{eqnarray*}
in $H^1(S\ot S/R,\dul{\Pic})$. Applying $H^{-1}$ to both sides, we find that
$$[\Can_R(S;S)\ot \Cc]=[(\Cc\ot_{S^{\ot 2}}\Ee)\ot\Can_R(S;S)]$$
in $\Br^c(S\ot S/R)$. Since the inverse of $[\Can_R(S;S)\ot \Cc]$ in $\Br^c(S\ot S/R)$
is $[\Can_R(S;S)\ot \Cc^*]$, we obtain that
$$[(\Can_R(S;S)\ot \Cc^*)\ot_{S^{\ot 4}} ((\Cc\ot_{S^{\ot 2}}\Ee)\ot\Can_R(S;S))]=
[(\Cc\ot_{S^{\ot 2}}\Ee)\ot \Cc^*]=1$$
in $\Br^c(S\ot S/R)$. Consequently $(\Cc\ot_{S^{\ot 2}}\Ee)\ot \Cc^*=\Ff$
is an elementary coring, and
$$(\Cc\ot_{S^{\ot 2}}\Ee)\ot \Cc^*\ot\Cc=\Ff\ot \Cc.$$
We have seen in \prref{5.2} that $\Cc\ot\Cc^*$ is elementary, and it follows that
$\Cc\ot_{S^{\ot 2}}\Ee\sim \Cc$.
\end{proof} 

\begin{lemma}\lelabel{5.5a}
Let $f:\ S\to T$ be a morphism of faithfully flat commutative $R$-algebras. If $\Cc$
is an Azumaya $S/R$-coring, then $\Cc\sim \tilde{f}(\Cc)=\Cc\ot_{S^{\ot 2}}
\Can_R(T;T)$.
\end{lemma}

\begin{proof}
As before, let $H(\Cc)=(\Cc,\alpha)$. Consider the maps $\varphi,\psi:\
S\to S\ot T$ given by
$$\varphi(s)=1\ot f(s)~~;~~\psi(s)=s\ot 1.$$
Applying \prref{1.2a.3}, we find that
\begin{eqnarray*}
&&\hspace*{-2cm}
[\varphi_*(\Cc,\alpha)]=[(S^{\ot 2}\ot (\Cc\ot_{S^{\ot 2}}T^{\ot 2}),\lambda_S\ot
(\alpha\ot_{S^{\ot 3}}T^{\ot 3}))]\\
&=& [\psi_*(\Cc,\alpha)]=[(\Cc\ot T^{\ot 2},\alpha\ot T^{\ot 3})]
\end{eqnarray*}
in $H^1(S\ot T/R,\dul{\Pic})$. Consequently
\begin{eqnarray*}
&&\hspace*{-2cm}
[H^{-1}(\varphi_*(\Cc,\alpha))]=[\Can_R(S;S)\ot (\Cc\ot_{S^{\ot 2}}\Can_R(T;T))]\\
&=& [H^{-1}(\psi_*(\Cc,\alpha))]=[\Cc\ot \Can_R(T;T)]
\end{eqnarray*}
in $\Br^c(S\ot T/R,\dul{\Pic})$. The inverse of $[\Cc\ot \Can_R(T;T)]$ in $\Br^c(S\ot T/R,\dul{\Pic})$
is $[\Cc^*\ot \Can_R(T;T)]$, and it follows that
\begin{eqnarray*}
&&\hspace*{-2cm}
(\Cc^*\ot \Can_R(T;T))\ot_{S\ot S\ot T\ot T}
(\Can_R(S;S)\ot (\Cc\ot_{S^{\ot 2}}\Can_R(T;T)))\\
&\cong &\Cc^*\ot (\Cc\ot_{S^{\ot 2}}\Can_R(T;T))
\cong \Ee
\end{eqnarray*}
where $\Ee$ is an elementary $S\ot T/R$-coring. We then have
$$\Cc\ot \Cc^*\ot (\Cc\ot_{S^{\ot 2}}\Can_R(T;T))\cong
\Cc\ot\Ee.$$
We know from \prref{5.2} that $\Cc\ot\Cc^*$ is elementary, so we can conclude that
$\Cc\sim \tilde{f}(\Cc)=\Cc\ot_{S^{\ot 2}}
\Can_R(T;T)$.
\end{proof}

\begin{proposition}\prlabel{5.6}
Let $S$ be a commutative faithfully flat $R$-algebra. We have a well-defined
group monomorphism
$$i_S:\ \Br^{\rm c}(S/R)\to \Br^{\rm c}_{\rm fl}(R),~~
i_S([\Cc])=[\Cc].$$
If $f:\ S\to T$ is a morphism of commutative faithfully flat $R$-algebras,
then we have a commutative diagram
$$\begin{diagram}
\Br^c(S/R)&\rTo^{i_S}&\Br^c_{\rm fl}(R)\\
\dTo_{\tilde{f}}&\NE^{i_T}&\\
\Br^c(T/R).&&
\end{diagram}$$
\end{proposition}

\begin{proof}
It follows from \leref{5.5} that $i_S$ is well-defined. Let us show that
$i_S$ is a group homomorphism. Consider two Azumaya $S/R$-corings $\Cc$ and
$\Dd$. Then the $S\ot S/R$-coring $\Cc^*\ot \Cc=\Ee_1$ and
the $S/R$-coring $\Cc\ot_{S^{\ot 2}} \Cc^*=\Ee_2$ are both elementary.
From \leref{5.5}, it follows that
\begin{eqnarray*}
\Cc\ot \Dd&\sim & (\Cc\ot \Dd)\ot_{S^{\ot 4}} (\Cc^*\ot \Cc)\\
&\cong & (\Cc\ot_{S^{\ot 2}} \Cc^*)\ot (\Dd\ot_{S^{\ot 2}} \Cc)\\
&\sim &\Dd\ot_{S^{\ot 2}} \Cc\cong \Cc\ot_{S^{\ot 2}} \Dd.
\end{eqnarray*}
Consequently
$$i_S[\Cc\ot_{S^{\ot 2}} \Dd]=[\Cc\ot \Dd]=i_S[\Cc]i_S[\Dd].$$
It is clear that $i_S$ is injective.\\
Finally, it follows from \leref{5.5a} that
 $i_S[\Cc]=[\Cc]=[\Cc\ot_{S^{\ot 2}} \Can_R(T;T)]=(i_T\circ \tilde{f})[\Cc]$.
\end{proof}

\begin{theorem}\thlabel{5.7}
Let $R$ be a commutative ring. Then
$$ \Br_{\rm fl}^c(R)\cong \colim \Br^c(\bullet/R)\cong H^2(R_{\rm fl},\units).$$
\end{theorem}

\begin{proof}
It follows from \prref{5.6} and the definition of colimit that we have a
map
$$i:\ \colim \Br^c(\bullet/R)\to \Br_{\rm fl}^c(R).$$
Suppose that $A$ is an abelian group, and suppose that we have a collection
of maps $\alpha_S:\ \Br^{\rm c}(S/R)\to A$ such that
$$\alpha_T\circ \tilde{f}=\alpha_S,$$
for every morphism of faithfully flat commutative $R$-algebras $f:\ S\to T$.
Take $x\in \Br^c_{\rm fl}(R)$. Then $x$ is represented by an Azumaya $S/R$-coring
$\Cc$. We claim that the map
$$\alpha:\ \Br^c_{\rm fl}(R)\to A,~~\alpha(x)=\alpha_S[\Cc]$$
is well-defined. Take an Azumaya $T/R$-coring $\Dd$ that also represents $x$.
Then 
$$\Cc\ot \Can_R(T;T)\sim \Cc\sim \Dd\sim \Dd\ot \Can_R(S;S)$$
and it follows from the injectivity of $i_{S\ot T}$ (see \prref{5.6}) that
$[\Cc\ot \Can_R(T;T)]=[\Dd\ot \Can_R(S;S)]$ in $\Br^c(S\ot T/R)$, hence
$$\alpha_S[\Cc]=\alpha_{S\ot T}[\Cc\ot \Can_R(T;T)]=\alpha_{S\ot T}
[\Dd\ot \Can_R(S;S)]=\alpha_T[\Dd],$$
as needed. We have constructed $\alpha$ in such a way that the
diagrams
$$\begin{diagram}
\Br^c(S/R)&\rTo^{i_S}&\Br^c_{\rm fl}(R)\\
&\SE_{\alpha_S}&\dTo_{\alpha}\\
&&A
\end{diagram}$$
commute.
This means that $\Br^c_{\rm fl}(R)$ satisfies the required universal property.
Finally, apply \equref{3.7.2}.
\end{proof}

\begin{corollary}\colabel{5.8}
Let $S$ be a faithfully flat commutative $R$-algebra. Then
$$\Ker(\Br_{\rm fl}^c(R)\to \Br_{\rm fl}^c(S))=
\Br^c(S/R).$$
\end{corollary}

\begin{proof}
Applying \coref{3.7}, \equref{villa2}, \equref{et1} with $q=1$ and
\thref{5.7}, we find that
\begin{eqnarray*}
\Br^c(S/R)&\cong&
H^1(S/R,\dul{\Pic})\cong H^1(S/R,C^1)\\
&\cong& \Ker(H^2(R_{\rm fl},\units)\to 
H^2(S_{\rm fl},\units))\\
&\cong&
\Ker(\Br_{\rm fl}^c(R)\to \Br_{\rm fl}^c(S)).
\end{eqnarray*}
\end{proof}

All our results remain valid if we replace the condition that $S$ is 
faithfully flat by the condition that $S$ is an \'etale covering, a
faithfully projective extension or a Zarisky covering of $R$
(see e.g. \cite{KO1} for precise definitions). It follows from
Artin's Refinement Theorem \cite{A1} that the (injective) map
$$\check H^2(R_{\rm et},\units)\to H^2(R_{\rm et},\units)$$
is an isomorphism. We will now present an algebraic interpretation
of $\check H^2(R_{\rm fl},\units)$ independent of Artin's Theorem.
Consider the subgroup $\Br^{\rm cnb}_{\rm fl}(R)$ consisting of classes
of Azumaya corings represented by an Azumaya coring with normal basis.

\begin{theorem}\thlabel{5.9}
Let $R$ be a commutative ring. Then
$$\Br^{\rm cnb}_{\rm fl}(R)\cong \check H^2(R_{\rm fl},\units).$$
\end{theorem}

\begin{proof}
Let $S$ be a faithfully flat commutative $R$-algebra, and consider the
maps
$$\begin{diagram}
&&\Br^c(S/R)\\
&&\dTo^{\gamma}\\
H^2(S/R,\units)&\rTo^{\beta}&H^1(S/R,\dul{\Pic})\\
\dTo{}&&\dTo^{}\\
\check H^2(R_{\rm fl},\units)&\hookrightarrow&H^2(R_{\rm fl},\units).
\end{diagram}$$
If $\Cc$ is an Azumaya $S/R$-coring with normal basis, then
$\gamma[\Cc]\in \im(\beta)$, so the image of $[\gamma]$ in 
$H^2(R_{\rm fl},\units)$ lies in the subgroup $\check H^2(R_{\rm fl},\units)$.
It follows that we have a monomorphism $\kappa:\ \Br^{\rm cnb}_{\rm fl}(R)
\to \check H^2(R_{\rm fl},\units)$ such that the following diagram
commutes:
$$\begin{diagram}
\Br^{\rm cnb}_{\rm fl}(R)&\hookrightarrow & \Br^{\rm c}_{\rm fl}(R)\\
\dTo^{\kappa}&&\dTo_{\cong}\\
\check H^2(R_{\rm fl},\units)&\hookrightarrow& H^2(R_{\rm fl},\units).
\end{diagram}$$
It is clear that $\gamma$ is surjective.
\end{proof}

\end{document}